\documentclass[11pt,fleqn]{article}

\usepackage{amsmath,amssymb,amsthm,enumerate}

\newcommand{\p}{\partial}
\newcommand{\const}{\mathop{\rm const}\nolimits}
\newcommand{\diag}{\mathop{\rm diag}\nolimits}

\newtheorem{theorem}{Theorem}
\newtheorem{lemma}{Lemma}
\newtheorem{corollary}{Corollary}
\newtheorem{proposition}{Proposition}
{\theoremstyle{definition} \newtheorem{definition}{Definition}
\newtheorem{example}{Example}

\newtheorem{note}{Note}

\begin{document}

\par\noindent {\LARGE\bf Reduction Operators \\of Linear Second-Order Parabolic Equations
\par}
{\it \vspace{4mm}\par\noindent 
Roman~O.~Popovych
\par\vspace{2mm}\par\vspace{2mm}\par\noindent 
Institute of Mathematics of National Academy of Sciences of Ukraine,
3 Tereshchenkivska Str., \\ Kyiv-4, 01601 Ukraine\\
Fakult\"at f\"ur Mathematik, Universit\"at Wien, Nordbergstra{\ss}e 15, A-1090 Wien, Austria\\
{\rm E-mail: }rop@imath.kiev.ua}




{\vspace{5mm}\par\noindent\hspace*{5mm}\parbox{150mm}{\small
The reduction operators, i.e., the operators of nonclassical (conditional) symmetry, 
of $(1+1)$-dimensional second order linear parabolic partial differential equations 
and all the possible reductions of these equations to ordinary differential ones 
are exhaustively described. 
This problem proves to be equivalent, in some sense, to solving the initial equations. 
The ``no-go'' result is extended to the investigation of point transformations 
(admissible transformations, equivalence transformations, Lie symmetries) and Lie reductions 
of the determining equations for the nonclassical symmetries. 
Transformations linearizing the determining equations are obtained in the general case and under 
different additional constraints. 
A nontrivial example illustrating applications of reduction operators to finding exact solutions 
of equations from the class under consideration is presented. 
An observed connection between reduction operators and Darboux transformations is discussed. 
}\par\vspace{3mm}}

\section{Introduction}

The notion of nonclassical symmetry (called also \mbox{$Q$-conditional} or, simply, conditional symmetry)
was introduced in~\cite{Bluman&Cole1969} by the example of the $(1+1)$-dimensional linear heat equation 
and a particular class of operators.
A precise and rigorous definition was suggested later 
(see, e.g.,~\cite{Fushchych&Tsyfra1987,Fushchych&Zhdanov1992,Zhdanov&Tsyfra&Popovych1999}).
In contrast to classical Lie symmetry, the system of determining equations on the coefficients 
of conditional symmetry operators of the heat equation was found to be nonlinear and less overdetermined~\cite{Bluman&Cole1969}.
First this system was investigated in~\cite{Webb1990} in detail, 
where it was partially linearized and its Lie symmetries were found. 
The problem on conditional symmetries of the heat equation was completely solved 
in \cite{Fushchych&Shtelen&Serov&Popovych1992}, see also~\cite{Fushchych&Shtelen&Serov1993en}. 
Namely, the determining equations were obtained in both the cases arising under consideration and then 
studied from the Lie symmetry point of view and reduced to the initial equation with nonlocal transformations.
The maximal Lie invariance algebras of the both sets of the determining equations appeared isomorphic to 
the maximal Lie invariance algebra of the initial equation. 
(Well later few of these results were re-obtained in~\cite{Mansfield1999}.)
The results of \cite{Fushchych&Shtelen&Serov&Popovych1992} were extended 
in \cite{Fushchych&Popowych1994-1,Popovych1995,Popovych1997} to a class 
of linear transfer equations which generalize the heat equation. 
Thus, for these equations the ``no-go'' theorems on linearization of determining equations for 
coefficients of conditional symmetry operators to the initial equations were proved in detail 
and wide multi-parametric families of exact solutions were constructed with non-Lie reductions. 
It was observed in \cite{Zhdanov&Lahno1998} that 
the proof of the theorem from \cite{Fushchych&Shtelen&Serov&Popovych1992} 
on reducibility of determining equations to initial ones in case of conditional symmetry operators 
with vanishing coefficients of~$\p_t$ are extended to the class of $(1+1)$-dimensional evolution equations. 
This theorem was also generalized to multi-dimensional evolution equations~\cite{Popovych1998} and 
even systems of such equations~\cite{Vasilenko&Popovych1999}.

The conditional invariance of a differential equation with respect to an involutive family of $l$~vector fields is equivalent to that 
any ansatz associated with this family reduces the equation to a differential equation with the lesser by $l$ number 
of independent variables~\cite{Zhdanov&Tsyfra&Popovych1999}. 
That is why, we use the shorter and more natural term ``reduction operators'' 
instead of ``operators of conditional symmetry'' or ``operators of nonclassical symmetry''
and say that a family of operators reduces a differential equation in case of
the equation is reduced by the associated ansatz.

In this paper we investigate the reduction operators of 
second-order linear parabolic partial differential equations in two independent variables, 
which have the general form
\begin{equation}\label{EqGenLPE}
Lu=u_t-A(t,x)u_{xx}-B(t,x)u_x-C(t,x)u=0,  
\end{equation}
where the coefficients $A$, $B$ and $C$ are (real) analytic functions of~$t$ and~$x$, $A\ne 0$. 
These coefficients form the entire tuple of arbitrary elements of class~\eqref{EqGenLPE}. 
We justify the partition of the sets of reduction operators into two subsets depending on 
vanishing or nonvanishing of the coefficients of $\p_t$. 
Usually this point is missed in the literature on conditional symmetries.
After factorization by the equivalence relation between reduction operators, we find 
the determining equations for the coefficients of operators from both the subsets. 
All the possible reductions of equations from class~\eqref{EqGenLPE} to ordinary differential equations are described. 
Different kinds of ``no-go'' statements on reduction of study (including solution) of the determining equations  
to the corresponding initial ones are obtained for equations from class~\eqref{EqGenLPE}. 
In particular, the point transformations of all kinds in both the classes of determining equations 
(admissible transformations, transformations from the associated equivalence groups, Lie symmetry transformations) 
are induced by the corresponding point transformations in class~\eqref{EqGenLPE}. 
Lie solutions of the determining equations prove to admit nontrivial interpretations 
in terms of Lie invariance properties of the initial equations.
An example on application of reduction operators is presented. 
It shows that in spite of the ``no-go'' statements nonclassical symmetry is 
an effective tool for finding exact solutions of partial differential equations. 

There are a number of motivations inducing us to carry out the above investigations. 
Class~\eqref{EqGenLPE} contains important subclasses that are widely applied in different science 
(probability theory, physics, financial mathematics, biology, etc.).
The most famous examples are the Kolmogorov equations ($C=0$) and adjoint to them the Fokker--Planck equations ($A_{xx}-B_x+C=0$)
which form a basis for analytical methods in the investigation of continuous-time continuous-state Markov processes. 
(The other names are Kolmogorov backward and Kolmogorov forward equations, respectively.)
The first use of the Fokker--Planck equation was the statistical description of Brownian motion of a particle in a fluid. 
Fokker--Planck equations with different coefficients 
also describe the evolution of one-particle distribution functions of a dilute gas with long-range collisions, 
problems of diffusion in colloids, population genetics, stock markets, quantum chaos, etc. 
Due to their importance and relative simplicity, equations from class~\eqref{EqGenLPE} are conventional objects 
for studies in the framework of group analysis of differential equations. 
Lie symmetries of these equations were classified by S. Lie~\cite{Lie1881}. 
The (1+1)-dimensional linear heat equation is often used as an illustrative example in textbooks on the subject \cite{Olver1986}
and a benchmark example for computer programs calculating symmetries of differential equations \cite{Head}. 
It is the equation that is connected with the invention of nonclassical symmetries \cite{Bluman&Cole1969}. 
First discussions on weak symmetries also involved the linear heat equation 
and a Fokker--Planck equation~\cite{Olver&Rosenau1987,Pucci&Saccomandi1992}. 
At the same time, all previous studies of nonclassical symmetries of equations~\eqref{EqGenLPE} were not systematic. 
Only a few equations and single properties were considered.  

Results of \cite{Fushchych&Popowych1994-1,Fushchych&Shtelen&Serov&Popovych1992,Popovych1995,Popovych2006b} 
are extended in the present paper mainly in two directions. 
Firstly, the entire class~\eqref{EqGenLPE} is regularly investigated with the nonclassical symmetry point of view and, 
secondly, non-evident properties of point transformations and Lie reductions of the determining equations are found 
via involving admissible transformations in the framework of nonclassical symmetries. 

Our paper is organized as follows: 
Necessary notions and statements on nonclassical symmetries are presented in Section~\ref{SectionOnDefOfRedOps}. 
The notion of equivalence of nonclassical symmetries with respect to a transformation group or a set of admissible transformations 
plays a crucial role in our consideration and therefore is separately given in Section~\ref{SectionOnEquivOfRedOpsWrtTransGroup}.
Section~\ref{SectionOnLieGroupAnalysisOfLSPEs} is devoted to reviewing 
known results on admissible transformations, point symmetries and equivalences in class~\eqref{EqGenLPE}, including discrete ones. 
The presentation of these results is important since they form a basis for the application 
of our technique involving transformations between equations 
and are extended in the paper to both the classes of determining equations. 
Moreover, Lie symmetry operators are special cases of reduction operators. 
The determining equations are derived in Section~\ref{SectionOnDetEqsForReductionOpsOfLPEs} 
for both the cases of nonvanishing and vanishing coefficients of $\p_t$. 
It is proved in Section~\ref{SectionOnLinearizationOfDetEqsForredictionOpsOfLPEs} via description of all possible reductions 
that solving the determining equations is equivalent to the construction of parametric families of solution
of the corresponding initial equations. 
As a result, nonlocal transformations reducing the determining equations to the initial ones are found. 
Point transformations and Lie reductions of the determining equations are studied in 
Sections~\ref{SectionOnLieSymsAndEquivGroupsOfDeterminingEqs} and~\ref{SectionOnLieReductionsOfDetEqsOfReditionOpsOfLPEs}, respectively.
Results on Lie reductions of the determining equations corresponding to reduction operators with zero coefficients of~$\p_t$ 
are presented in such a form that they are directly extended to the general class of $(1+1)$-dimensional evolution equations. 
In Section~\ref{SectionOnParticularCasesOfReductionsAndLinearization} we investigate the determining equations 
along with some non-Lie additional constraints. 
A nontrivial application of reduction operators to finding exact solutions of equations from class~\eqref{EqGenLPE}, arising under 
Lie reductions of the Navier--Stokes equations, is presented in Section~\ref{SectionOnReductionOpsOfLinearTransferEqs}.
In the last section we discuss possible extensions of obtained results, in particular, via study of the observed connection between 
reduction operators and the Darboux transformations of equations from class~\eqref{EqGenLPE}. 

To check  results on Lie invariance of differential equations appearing in the paper, 
we used the unique program LIE by A.~Head~\cite{Head}.

\section{Reduction operators of differential equations}
\label{SectionOnDefOfRedOps}

Following~\cite{Fushchych&Tsyfra1987,Fushchych&Zhdanov1992,Popovych&Vaneeva&Ivanova2007,Zhdanov&Tsyfra&Popovych1999}, 
in this section we shortly adduced necessary notions and results on nonclassical (conditional) symmetries of 
differential equations. 
After substantiating with different arguments, we use the name ``families of reduction operators'' 
instead of ``involutive families of nonclassical (conditional) symmetry operators''. 
  
Consider an involutive family $Q=\{Q^1,\ldots,Q^l\}$ of $l$ ($l\leqslant n$) first order differential operators 
\[
Q^s=\xi^{si}(x,u)\p_i+\eta^s(x,u)\p_u, \quad s=1,\dots,l,
\]
in the space of the variables~$x$ and~$u$, satisfying the condition $\mathop{\rm rank}\nolimits \|\xi^{si}(x,u)\|=l$.

Hereafter 
$x$ denote the $n$-tuple of independent variables $(x_1,\ldots,x_n)$ 
and $u$ is treated as the unknown function. 
The index $i$ runs from 1 to $n$,  
the indices $s$ and $\sigma$ run from 1 to $l$, 
and we use the summation convention for repeated indices.
$\p_i=\p/\p x_i$, $\p_u=\p/\p u$.
Subscripts of functions denote differentiation with respect to the corresponding variables. 
The local consideration is assumed.

The requirement of involution for the family~$Q$ means that the commutator of any pair of operators from~$Q$ 
belongs to the span of $Q$ over the ring  of smooth functions of the variables~$x$ and $u$, i.e., 
\[
\forall\, s,s'\quad \exists \, \zeta^{ss'\sigma}=\zeta^{ss'\sigma}(x,u)\colon\quad [Q^s,Q^{s'}]=\zeta^{ss'\sigma}Q^\sigma.\]
The set of such families will be denoted by $\mathcal Q^l$.

If operators~$Q^1$, \ldots, $Q^l$ form an involutive family~$Q$, then the family $\widetilde Q$ 
of differential operators
\[
\widetilde Q^s=\lambda^{s\sigma}Q^\sigma, 
\qquad\mbox{where}\qquad \lambda^{s\sigma}=\lambda^{s\sigma}(x,u),\quad \det\|\lambda^{s\sigma}\|\not=0,
\]
is also involutive and is called {\em equivalent} to the family~$Q$. 
This will be denoted by $\widetilde Q=\{\widetilde Q^s\}\sim Q=\{Q^s\}$. 
(In the case $l=1$  the functional matrix~$(\lambda^{s\sigma})$ becomes a single nonvanishing multiplier $\lambda=\lambda(x,u)$.)
Denote also the result of factorization of~$\mathcal Q^l$ with respect to this equivalence relation by~$\mathcal Q^l_{\rm f}$. 
Elements of~$\mathcal Q^l_{\rm f}$ will be identified with their representatives in~$\mathcal Q^l$. 

If a family consists of a single operator ($l=1$), the involution condition degenerates to an identity. 
Therefore, in this case we can omit the words ``involutive family'' and talk only about operators. 
Thus, two differential operators are equivalent if they differ on a multiplier 
being a non-vanishing function of~$x$ and~$u$.

The first order differential function~$Q^s[u]:=\eta^{s}(x,u)-\xi^{si}(x,u)u_i$ 
is called the {\it characteristic} of the operator~$Q^s$.
In view of the Frobenius theorem, the above involution condition is equivalent to 
that the characteristic system $Q[u]=0$ of PDEs $Q^s[u]=0$ 
(called also the \emph{invariant surface condition}) has 
$n+1-l$ functionally independent integrals $\omega^0(x,u),\ldots,\omega^{n-l}(x,u)$. 
Therefore, the general solution of this system can be implicitly presented in the form 
$F(\omega^0,\ldots,\omega^{n-l})=0$, where~$F$ is an arbitrary function of its arguments.

The characteristic systems of equivalent families of operators have the same set of solutions.
And vice versa, any family of $n+1-l$ functionally independent functions of~$x$ and $u$ 
is a complete set of integrals of 
the characteristic system of an involutive family of $l$~differential operators. 
Therefore, there exists the one-to-one correspondence between $\mathcal Q^l_{\rm f}$ 
and the set of families of $n+1-l$ functionally independent functions of~$x$ and $u$, 
which is factorized with respect to the corresponding equivalence. 
(Two families of the same number of functionally independent functions of the same arguments are considered equivalent 
if any function from one of the families is functionally dependent on functions from the other family.)

A function $u=f(x)$ is called \emph{invariant with respect to the involutive operator family}~$Q$ 
(or, briefly, $Q$-\emph{invariant}) if 
it is a solution of the characteristic system $Q[u]=0$. 
This notion is justified by the following facts. 
Any involutive family of $l$ operators is equivalent to a basis $\widetilde Q=\{\widetilde Q^s\}$ 
of an $l$-dimensional (Abelian) Lie algebra $\mathfrak g$ of vector fields in the space $(x,u)$. 
Each solution $u=f(x)$ of the associated characteristic system satisfies the characteristic system $\widetilde Q[u]=0$. 
Therefore, the graph of the function $u=f(x)$ is invariant with respect to the $l$-parametric local transformation group 
generated by the algebra~$\mathfrak g$.

Since $\mathop{\rm rank}\nolimits \|\xi^{si}(x,u)\|=l$, we can assume without loss of generality that 
$\omega^0_u\not=0$ and $F_{\omega^0}\not=0$ and resolve the equation~$F=0$ with respect to~$\omega^0$:
$
\omega^0=\varphi(\omega^1,\ldots,\omega^{n-l})$.
This representation of the function~$u$ is called an \emph{ansatz} corresponding to the family~$Q$.

Consider an $r$th-order differential equation~$\mathcal L$ of the form~$L(x,u_{(r)})=0$
for the unknown function $u$ of $n$ independent variables $x=(x_1,\ldots,x_n).$
Here, $u_{(r)}$ denotes the set of all the derivatives of the function $u$ with respect to $x$
of order not greater than~$r$, including $u$ as the derivative of the zero order.
Within the local approach the equation~$\mathcal L$ is treated as an algebraic equation 
in the jet space $J^{(r)}$ of the order $r$ and is identified with the manifold of its solutions in~$J^{(r)}$. 
Denote this manifold by the same symbol~$\mathcal L$ and 
the manifold defined by the set of all the differential consequences of the characteristic system~$Q[u]=0$ 
in $J^{(r)}$ by $\mathcal Q^{(r)}$, i.e., 
\begin{gather*}
\mathcal Q^{(r)}=\{ (x,u_{(r)}) \in J^{(r)}\, |\, D_1^{\alpha_1}\ldots D_n^{\alpha_n}Q^s[u]=0, 
\;\: \alpha_i\in\mathbb{N}\cup\{0\},\;\:|\alpha|\mbox{:}=\alpha_1+\cdots+\alpha_n<r \},
\end{gather*}
where $D_i=\p_{x_i}+u_{\alpha+\delta_i}\p_{u_\alpha}$ is the operator of total differentiation with respect to the variable~$x_i$, 
$\alpha=(\alpha_1,\ldots,\alpha_n)$ is an arbitrary multi-index,  
$\delta_i$ is the multiindex whose $i$th entry equals 1 and whose other entries are zero. 
The variable $u_\alpha$ of the jet space $J^{(r)}$ corresponds to the derivative
$\p^{|\alpha|}u/\p x_1^{\alpha_1}\ldots\p x_n^{\alpha_n}$.

\begin{definition}\label{DefinitionOfCondSym}
The differential equation~$\mathcal L$ is called \emph{conditionally invariant} with respect to 
the involutive family $Q$ if the relation $Q^s_{(r)}L(x,u_{(r)})\bigl|_{\mathcal L\cap\mathcal Q^{(r)}}=0$ holds,
which is called the \emph{conditional invariance criterion}.
Then $Q$ is called an \emph{involutive family of conditional symmetry} 
(or $Q$-conditional symmetry, nonclassical symmetry etc) operators of the equation~$\mathcal L$.
Here the symbol $Q^s_{(r)}$ stands for the standard $r$th prolongation
of the operator~$Q^s$ \cite{Olver1986,Ovsiannikov1982}:
$Q^s_{(r)}=Q^s+\sum_{|\alpha|{}\leqslant  r} \eta^{s\alpha}\p_{u_\alpha}$, 
where 
$\eta^{s\alpha}=D_1^{\alpha_1}\ldots D_n^{\alpha_n}Q^s[u]+\xi^{si}u_{\alpha+\delta_i}$.
\end{definition}

The equation~$\mathcal L$ is conditionally invariant with respect to the family~$Q$
if and only if the ansatz constructed with this family reduces~$\mathcal L$
to a differential equation with $n-l$ independent variables~\cite{Zhdanov&Tsyfra&Popovych1999}.
So, we will also call involutive families of conditional symmetry operators 
{\it families of reduction operators} of~$\mathcal L$. 
Another treatment of conditional invariance is that the system $\mathcal L\cap\mathcal Q^{(r)}$ is compatible 
in the sense of absence of nontrivial differential consequences~\cite{Olver1994,Olver&Vorob'ev1996}. 
Note that the paper~\cite{Olver1994} contains a number of interesting statements and ideas on the subject, 
which, unfortunately, did not become known according to their merits. 
If the infinitesimal invariance condition is not satisfied but nevertheless 
the equation~$\mathcal L$ has $\mathcal Q$-invariant solutions 
then $\mathcal Q$ is called a family of weak symmetry operators of the equation~$\mathcal L$ 
\cite{Olver&Rosenau1987,Olver&Vorob'ev1996}. 
Nonclassical symmetries are often defined as generators of parametric groups of transformations 
preserving the solutions of~$\mathcal L$ which additionally satisfy the corresponding invariant surface condition \cite{Hydon2000}. 
It is necessary to precisely interpret all the terms involved in this definition since otherwise it leads to the conclusion that, roughly speaking,   
any operator is a nonclassical symmetry of any partial differential equation. 
See also \cite{Bila&Niesen2004,Clarkson1995,Olver&Vorob'ev1996} for the discussion of connections between 
different kinds of symmetries. 

\begin{lemma}[\cite{Fushchych&Zhdanov1992,Zhdanov&Tsyfra&Popovych1999}]\label{LemmaOnEquivFamiliesOfOperators} 
If a differential equation  is conditionally invariant with respect to an operator family~$Q$, 
then it is conditionally invariant with respect to any family of operators, which is equivalent to~$Q$. 
\end{lemma}

The set of involutive families of $l$ reduction operators of the equation~$\mathcal L$ 
is a subset of $\mathcal Q^l$ and so will be denoted by $\mathcal Q^l(\mathcal L)$. 
In view of Lemma~\ref{LemmaOnEquivFamiliesOfOperators}, 
$Q\in \mathcal Q^l(\mathcal L)$ and $\widetilde Q\sim Q$ imply $\widetilde Q\in \mathcal Q^l(\mathcal L)$, 
i.e., $\mathcal Q^l(\mathcal L)$ is closed under the equivalence relation on $\mathcal Q^l$.
Therefore, factorization of $\mathcal Q^l$ with respect to this equivalence relation can be naturally restricted  
on~$\mathcal Q^l(\mathcal L)$ that results in the subset~$\mathcal Q^l_{\rm f}(\mathcal L)$ of $\mathcal Q^l_{\rm f}$. 
As in the whole set~$\mathcal Q^l_{\rm f}$, we identify elements of~$\mathcal Q^l_{\rm f}(\mathcal L)$ 
with their representatives in~$\mathcal Q^l(\mathcal L)$.
In this approach the problem of complete description of families of $l$ reduction operators 
for the equation~$\mathcal L$ is nothing but the problem of finding~$\mathcal Q^l_{\rm f}(\mathcal L)$. 

A different terminology can be used to call elements of~$\mathcal Q^l_{\rm f}$. 
Namely, it is possible to consider each element of~$\mathcal Q^l_{\rm f}$ as a $C^\infty$-module of the module dimension~$l$, 
closed with respect to commutation~\cite{Olver&Vorob'ev1996,Vorob'ev1991}.

There are families of reduction operators related to classical Lie symmetries. 
Let $\mathfrak g$ be an $l$-dimensional Lie invariance algebra of the equation~$\mathcal L$, 
whose basis operators satisfy the condition
$\mathop{\rm rank}\nolimits \|\xi^{si}\|=  
\mathop{\rm rank}\nolimits \|\xi^{si},\eta^s\|\ ({}=l'\leqslant l).$
The subsets consisting of $l'$ elements of~$\mathfrak g$, 
which are linearly independent over the ring of smooth functions of $x$ and $u$, belong to~$\mathcal Q^{l'}(\mathcal L)$ 
and are equivalent each to other. 
The families of similar kind and ones equivalent to them will be called \emph{Lie families of reduction operators}. 
The other families of reduction operators will be called \emph{non-Lie}.

\section{Equivalence of families of reduction operators\\ with respect to transformation groups}
\label{SectionOnEquivOfRedOpsWrtTransGroup}

\looseness=1
We can essentially simplify and order the investigation of reduction operators, additionally taking into account 
Lie symmetry transformations in case of a single equation~\cite{Popovych2000} and  
transformations from the equivalence group or the whole set of admissible transformations 
in case of a class of equations~\cite{Popovych&Vaneeva&Ivanova2007}. 
Then the problem becomes similar to group classification of differential equations.

\begin{lemma}
Any point transformation of $x$ and $u$ induces a one-to-one mapping of~$\mathcal Q^l$ into itself.
Namely, the transformation~$g$: $\tilde x=X(x,u)$, $\tilde u=U(x,u)$ generates 
the mapping~\mbox{$g_*^l\colon \mathcal Q^l\to\mathcal Q^l$} such that 
the involutive family~$Q$ is mapped to the involutive family $g_*^lQ$ consisting from the operators
$g_*Q^s=\tilde\xi^{si}\p_{\tilde x_i}+\tilde\eta^s\p_{\tilde u}$, where 
$\tilde\xi^{si}(\tilde x,\tilde u)=Q^sX^i(x,u)$, 
$\tilde\eta^s(\tilde x,\tilde u)=Q^sU(x,u)$. 
If~$Q'\sim Q$ then  $g_*^l Q'\sim g_*^lQ$. 
Therefore, the corresponding factorized mapping~$g_{\rm f}^l\colon\mathcal Q^l_{\rm f}\to\mathcal Q^l_{\rm f}$ also 
is well-defined and one-to-one.

\end{lemma}

\begin{definition}[\cite{Popovych&Korneva1998,Popovych2000}]\label{DefinitionOfEquivInvFamiliesWrtGroup}
Involutive families $Q$ and $\widetilde Q$ of the same number~$l$ of differential operators are called 
\emph{equivalent with respect to a group $G$ of point transformations} if there exists a transformation $g$ from $G$ 
for which the families $Q$ and $g_*^l\widetilde Q$ are equivalent.

\noindent
{\it Notation:} $Q\sim \widetilde Q \bmod G.$
\end{definition}

\begin{lemma}\label{LemmaOnInducedMappingOfRedictionOps}
Given any point transformation $g$ of the equation~$\mathcal L$ to an equation~$\tilde{\mathcal L}$,
$g_*^l$ maps~$\mathcal Q^l(\mathcal L)$ to~$\mathcal Q^l(\tilde{\mathcal L})$ in a one-to-one manner. 
The same statement is true for the factorized mapping $g_{\rm f}^l$ from $\mathcal Q^l_{\rm f}(\mathcal L)$
to~$\mathcal Q^l_{\rm f}(\tilde{\mathcal L})$.
\end{lemma}

\begin{corollary}\label{CorollaryOnEquivReductionOperatorWrtSymGroup}
Let $G$ be a Lie symmetry group of the equation~$\mathcal L.$ Then the equivalence of involutive families of 
$l$ differential operators with respect to the group $G$ generates equivalence relations in~$\mathcal Q^l(\mathcal L)$ 
and in~$\mathcal Q^l_{\rm f}(\mathcal L)$.  
\end{corollary}

Consider a class~$\mathcal L|_{\mathcal S}$ of equations~$\mathcal L_\theta$: 
$L(x,u_{(r)},\theta(x,u_{(r)}))=0$ parameterized by~$\theta$.
Here, $L$ is a fixed function of $x,$ $u_{(r)}$ and $\theta.$
The symbol~$\theta$ denotes the tuple of arbitrary (parametric) functions
$\theta(x,u_{(r)})=(\theta^1(x,u_{(r)}),\ldots,\theta^k(x,u_{(r)}))$
running through the solution set~$\mathcal S$ the system~$S(x,u_{(r)},\theta_{(q)}(x,u_{(r)}))=0$.
This system consists of differential equations on $\theta$,
where $x$ and $u_{(r)}$ play the role of independent variables
and $\theta_{(q)}$ stands for the set of all the partial derivatives of $\theta$ of order not greater than $q$.
In what follows we call the functions $\theta$ \emph{arbitrary elements}.
By $G^\sim$ we denote the point transformations group preserving the
form of the equations from~$\mathcal L|_{\mathcal S}$.

For a fixed value $l\leqslant n$, consider the set~$P=P(L,S)$ of all pairs each of which consists of
an equation $\mathcal L_\theta$ from~$\mathcal L|_{\mathcal S}$ and a family~$Q$ from~$\mathcal Q^l(\mathcal L_\theta)$.
In view of Lemma~\ref{LemmaOnInducedMappingOfRedictionOps},
the action of transformations from~$G^\sim$ on $\mathcal L|_{\mathcal S}$ and
$\{\mathcal Q^l(\mathcal L_{\theta})\,|\,\theta\in{\mathcal S}\}$
together with the pure equivalence relation of involutive families of $l$ differential operators
naturally generates an equivalence relation on~$P$.

\begin{definition}\label{DefinitionOfEquivOfRedOperatorsWrtEquivGroup}
Let $\theta,\theta'\in{\mathcal S}$,
$Q\in\mathcal Q^l(\mathcal L_\theta)$, $Q'\in\mathcal Q^l(\mathcal L_{\theta'})$.
The pairs~$(\mathcal L_\theta,Q)$ and~$(\mathcal L_{\theta'},Q')$
are called {\em $G^\sim$-equivalent} if there exists a transformation $g\in G^\sim$
which maps the equation~$\mathcal L_\theta$ to the equation~$\mathcal L_{\theta'}$, and
$Q'\sim g_*^lQ$.
\end{definition}

Classification of families of reduction operators with respect to~$G^\sim$ will be understood as
classification in~$P$ with respect to the above equivalence relation.
This problem can be investigated in a way similar to the usual group classification in classes
of differential equations.
Namely, we construct first the reduction operators which are defined for all values of the arbitrary elements.
Then we classify, with respect to the equivalence group, the values of arbitrary elements for which 
the corresponding equations admit additional families of reduction operators.

In an analogous way we also can introduce equivalence relations on~$P$, which are
generated by either generalizations of usual equivalence groups or
all admissible point transformations~\cite{Popovych&Eshraghi2005} 
(called also form-preserving ones~\cite{Kingston&Sophocleous1998})
in pairs of equations from~$\mathcal L|_{\mathcal S}$.

\begin{note}
The consideration of the previous and this sections and known examples of studying reduction operators  
lead to the empiric conclusion that 
possessing a wide Lie symmetry group by a differential equation~$\mathcal L$ complicates, in some way,  
finding nonclassical symmetries of~$\mathcal L$. 
Indeed, any subalgebra of the corresponding maximal Lie invariance algebra, satisfying the transversality condition, 
generates a class of equivalent Lie families of reduction operators. 
A non-Lie family of reduction operators existing, the action of symmetry transformations on it results in  
a series of non-Lie families of reduction operators, which are inequivalent in the usual sense. 
Therefore, for any fixed value of~$l$ the system of determining equations on coefficients of operators 
from $\mathcal Q^l(\mathcal L)$ is not sufficiently overdetermined to be completely integrated in an easy way, 
even after factorized with respect to the equivalence relation in $\mathcal Q^l(\mathcal L)$. 
To produce essentially different non-Lie reductions, one have to exclude the solutions of determining equations, 
which give Lie families of reduction operators and non-Lie families being equivalent to others with respect to 
the Lie symmetry group of~$\mathcal L$. 
As a result, the ratio of efficiency of such reductions to expended efforts can be vanishingly small.
\end{note}

\section{Lie group analysis of linear second-order parabolic equations}
\label{SectionOnLieGroupAnalysisOfLSPEs}

Group classification in class~\eqref{EqGenLPE} was first performed by S.~Lie~\cite{Lie1881} as a part of
his classification of general linear second-order PDEs in two independent variables.
(See also a modern treatment of this subject in~\cite{Ovsiannikov1982}.) 
We shortly adduce these classical results, extending them for our purposes with 
using the notions of admissible transformations and normalized classes of differential equations. 
First, normalization properties of different classes of linear second-order parabolic equations 
were simultaneously analyzed in~\cite{Popovych&Kunzinger&Ivanova2007} in detail. 

Roughly speaking, an admissible transformation in a class of systems of differential equations is 
a point transformation connecting at least two systems from this class (in the
sense that one system is transformed into the other by the transformation).
The equivalence group of the class is the set of admissible transformations which can be applied to every system from the class.
The class is called \emph{normalized} if any
admissible transformation in this class belongs to its equivalence group
and is called \emph{strongly normalized} if additionally the equivalence group
is generated by transformations from the point symmetry groups of systems from the class.
The set of admissible transformations of a \emph{semi-normalized class} is generated by
the transformations from the equivalence group of the whole class and the transformations
from the point symmetry groups of initial or transformed systems. 
Strong semi-normalization is defined in the same way as strong normalization.
Any normalized class is semi-normalized. 
Two systems from a semi-normalized class are transformed into one another by a point transformation
iff they are equivalent with respect to~the equivalence group of this class.
See \cite{Popovych2006a,Popovych2006c,Popovych&Eshraghi2005,Popovych&Kunzinger&Eshraghi2006} 
for precise definitions and statements.

Any point transformation~$\mathcal T$ in the space of variables $(t,x,u)$ has the form 
\[
\tilde t=\mathcal T^t(t,x,u),\quad 
\tilde x=\mathcal T^x(t,x,u),\quad 
\tilde u=\mathcal T^u(t,x,u), 
\]
where the Jacobian $|\p(\mathcal T^t,\mathcal T^x,\mathcal T^u)/\p(t,x,u)|$ does not vanish.

\begin{lemma}\label{LemmaOnAdmTransOfInhomLPEs}
A point transformation~$\mathcal T$ connects two equations from class~\eqref{EqGenLPE} if and only if  
$\,\mathcal T^t_x=\mathcal T^t_u=0$, $\mathcal T^x_u=0$, $\mathcal T^u_{uu}=0$, i.e.,
\begin{equation}\label{EqGenFormOfTransOfInhomLPEs}
\tilde t=T(t), \quad \tilde x=X(t,x), \quad \tilde u=U^1(t,x)u+U^0(t,x),
\end{equation}
where $T$, $X$, $U^1$ and $U^0$ are arbitrary smooth functions of their arguments such that \mbox{$T_tX_xU^1\ne0$} 
and additionally $U^0/U^1$ is a solution of the initial equation. 
The arbitrary elements are transformed by the formulas
\begin{gather}\label{EqTransOfCoeffsOfLPE}
\tilde A=\frac{X_x^2}{T_t}A,\quad 
\tilde B=\frac{X_x}{T_t}\left(B-2\frac{U^1_x}{U^1}A\right)-\frac{X_t-AX_{xx}}{T_t},\quad 
\tilde C=-\frac{U^1}{T_t}L\frac1{U^1}.
\end{gather}
Here $L=\p_t-A\p_{xx}-B\p_x-C$ is the second-order linear differential operator associated with 
the initial (non-tilde) equation.
\end{lemma}

\begin{corollary}\label{CorollaryOnSemiNormOfLPEs}
Class~\eqref{EqGenLPE} is strongly semi-normalized. 
The equivalence group~$G^\sim$ of class~\eqref{EqGenLPE} is formed by the transformations 
determined in the space of variables and arbitrary elements by formulas~\eqref{EqGenFormOfTransOfInhomLPEs} and  
\eqref{EqTransOfCoeffsOfLPE}, where $T$, $X$ and $U^1$ are arbitrary smooth functions of their arguments 
such that $\,T_tX_xU^1\ne0$ and $U^0=0$ additionally. 
\end{corollary}

\begin{note}
Due to the presence of the linear superposition principle, 
the class~\eqref{EqGenLPE} is not normalized because it is formed by linear homogeneous equations. 
The minimal normalized superclass of class~\eqref{EqGenLPE} is 
the associated class of inhomogeneous equations of the general form 
\[
u_t=A(t,x)u_{xx}+B(t,x)u_x+C(t,x)u+D(t,x).
\]
\end{note}

Using transformations from~$G^\sim$, the arbitrary elements~$A$ and~$B$ can be simultaneously gauged to~1 and~0, respectively. 
Hence, any equation from class~\eqref{EqGenLPE} can be 
reduced by a transformation from~$G^\sim$ to an equation of the general form
\begin{equation}\label{EqReducedLPE}
u_t-u_{xx}+V(t,x)u=0.
\end{equation}
The admissible transformations in the subclass~\eqref{EqReducedLPE} are those admissible
transformations in the class~\eqref{EqGenLPE} which preserve the gauges $A=1$ and $B=0$, 
i.e., which additionally satisfy the conditions $\mathcal T^t_t=(\mathcal T^x_x)^2$ and 
$2\mathcal T^x_x\mathcal T^u_{xu}=-\mathcal T^x_t\mathcal T^u_u$.

\begin{corollary}\label{CorollaryOnAdmTransOfLPEsA1B0}
A point transformation~$\mathcal T$ connects two equations from class~\eqref{EqReducedLPE} if and only if it has the form 
\begin{equation}\label{EqTransFromEquivGroupOfReducedLPEs}
\hspace*{-1\arraycolsep}\begin{array}{l}
\tilde t=\int\! \sigma^2dt,\quad \tilde x=\sigma x+\zeta,\quad 
\tilde u=U^1u+U^0,\quad
U^1:=\theta\exp\left(-\dfrac{\sigma_t}{4\sigma}x^2-\dfrac{\zeta_t}{2\sigma}x\right),
\\[3ex] 
\tilde V=\dfrac1{\sigma^2}\left(V+\dfrac{\sigma\sigma_{tt}-2\sigma_t{}^2}{4\sigma^2}x^2
+\dfrac{\sigma\zeta_{tt}-2\sigma_t\zeta_t}{2\sigma^2}x-\dfrac{\theta_t}\theta
-\dfrac{\sigma_t}{2\sigma}-\dfrac{\zeta_t{}^2}{4\sigma^2}\right),
\end{array}
\end{equation}
where $\sigma=\sigma(t)$, $\zeta=\zeta(t)$, $\theta=\theta(t)$ and $U^0=U^0(t,x)$ 
are arbitrary smooth functions of their arguments such that  
$\sigma\theta\ne0$ and $U^0/U^1$ is a solution of the initial equation.
Class~\eqref{EqReducedLPE} is strongly semi-normalized. 
Any transformation from the equivalence group~$G^\sim_{\rm r}$ of class~\eqref{EqReducedLPE} 
has the form~\eqref{EqTransFromEquivGroupOfReducedLPEs}, where $U^0=0$ additionally.
\end{corollary}

The narrower equivalence group under preserving certain normalization properties suggests
class~\eqref{EqReducedLPE} as the most convenient one for group classification. 
Moreover, solving the group classification problem for class~\eqref{EqGenLPE} is reduced 
to solving the group classification problem for class~\eqref{EqReducedLPE}.
The results on the group classification of class~\eqref{EqGenLPE} (resp.~\eqref{EqReducedLPE}) 
can be formulated in the form of the following theorem~\cite{Lie1881,Ovsiannikov1982}. 

\begin{theorem}\label{TheoremOnGroupClassificationOfLPEs}
The kernel Lie algebra of class~\eqref{EqGenLPE} (resp.~\eqref{EqReducedLPE}) is $\langle u\p_u\rangle$.
Any equation from class~\eqref{EqGenLPE} (resp.~\eqref{EqReducedLPE}) is invariant with respect 
to the operators $f\p_u$,
where the parameter-function $f=f(t,x)$ runs through the solution set of this equation.
All possible $G^\sim$-inequivalent (resp. $G^\sim_{\rm r}$-inequivalent) cases of extension of the maximal 
Lie invariance algebra are exhausted by the following ones 
(the values of~$V$ are given together with the corresponding maximal Lie invariance algebras):

\vspace{1.5ex}

$\makebox[6mm][l]{\rm 1.}
V=V(x) \colon\quad \langle\p_t,\,u\p_u,\,f\p_u\rangle;$

\vspace{1.5ex}

$\makebox[6mm][l]{\rm 2.}
V=\mu x^{-2},\ \mu\ne0 \colon\quad \langle\p_t,\, D,\, \Pi,\, u\p_u,\, f\p_u\rangle;$

\vspace{1.5ex}

$\makebox[6mm][l]{\rm 3.}
V=0 \colon\quad \langle\p_t,\, \p_x,\, G,\, D,\, \Pi,\, u\p_u,\, f\p_u\rangle$.

\vspace{1.5ex}

\noindent
Here $D=2t\p_t+x\p_x,\ \Pi=4t^2\p_t+4tx\p_x-(x^2+2t)u\p_u,\ G=2t\p_x-xu\p_u.$
\end{theorem}

Let $\mathcal L$ be an equation from class~\eqref{EqGenLPE}, 
$\mathfrak g(\mathcal L)$ denote its maximal Lie invariance algebra and 
$\mathfrak g^\infty(\mathcal L)$ be the infinite-dimensional ideal of this algebra, 
consisting of the operators of the form $f\p_u$, 
where the parameter-function $f=f(t,x)$ runs through the solution set of~$\mathcal L$.
The quotient algebra $\mathfrak g(\mathcal L)/\mathfrak g^\infty(\mathcal L)$ is identified with 
the finite-dimensional subalgebra $\mathfrak g^{\rm ess}(\mathcal L)$ of~$\mathfrak g(\mathcal L)$, 
spanned by the `essential' Lie invariance operators of~$\mathcal L$, which do not contain summands of the form $f(t,x)\p_u$.
Each operator from $\mathfrak g(\mathcal L)$ is similar to an operator from $\mathfrak g^{\rm ess}(\mathcal L)$ under 
a trivial linear-superposition transformation $\tilde t=t$, $\tilde x=x$, $\tilde u=u+f(t,x)$.

\begin{corollary}\label{CorollaryOnNumberOfNontrivSymsOfLPEs}
For every equation~$\mathcal L$ from class~\eqref{EqGenLPE} $\dim\mathfrak g^{\rm ess}(\mathcal L)\in\{1,2,4,6\}$. 
\end{corollary}

It will be shown below that for every equation~$\mathcal L$ from class~\eqref{EqGenLPE}
the number of reduction operators being inequivalent with respect to the group of linear-superposition transformations,
roughly speaking,  
is significantly greater that the number of `essential' Lie invariance operators.

\section{Determining equations for reduction operators\\ of linear second-order parabolic equations}
\label{SectionOnDetEqsForReductionOpsOfLPEs}

In the case of two independent variables $t$ and $x$ and one dependent variable $u$, each reduction operator is written as  
$Q=\tau(t,x,u)\p_t+\xi(t,x,u)\p_x+\eta(t,x,u)\p_u$, where $(\tau,\xi)\not=(0,0)$.
The conditional invariance criterion for an equation~$\mathcal L$ from class~\eqref{EqGenLPE} 
and the operator~$Q$ has the form~\cite{Fushchych&Tsyfra1987} 
\[
Q_{(2)}Lu\:\big|_{\;Lu=0,\;\; Q[u]=0,\;\; D_tQ[u]=0,\;\; D_xQ[u]=0}=0,
\]
where $Q_{(2)}$ is the standard second prolongation of~$Q$, $Q[u]=\eta-\tau u_t-\xi u_x$ is the characteristic of~$Q$, 
$D_t$ and $D_x$ denote the total differentiation operators with respect to~$t$ and~$x$:
\begin{gather*}
D_t=\p_t+u_t\p_u+u_{tt}\p_{u_t}+u_{tx}\p_{u_x}+\cdots, \\  
D_x=\p_x+u_x\p_u+u_{tx}\p_{u_t}+u_{xx}\p_{u_x}+\cdots.
\end{gather*}
All equalities hold true as algebraic relations in the second-order jet space~$J^{(2)}$ over the space of 
the independent variables $(t,x)$ and the dependent variable~$u$.

Since $\mathcal L$ is an evolution equation, 
there are two principally different cases of finding its reduction operators: $\tau\ne0$ and $\tau=0$. 
The investigation of these cases results in the preliminary description of the reduction operators.

\begin{lemma}\label{LemmaOnReductionOpsWithNonzeroCoeffsOfDtOfLPEs}
Every reduction operator of an equation~$\mathcal L$ from class~\eqref{EqGenLPE} is equivalent to either an operator 
\[\p_t+g^1(t,x)\p_x+(g^2(t,x)u+g^3(t,x))\p_u,\]
where the functions $g^1=g^1(t,x)$, $g^2=g^2(t,x)$ and $g^3=g^3(t,x)$ satisfy the system 
\begin{equation}\label{EqDet1ForRedOpsOfLPEs}\arraycolsep=0em
\begin{array}{l}
g^1_t-Ag^1_{xx}-Bg^1_x+\biggl(2g^1_x-\dfrac{A_x}{A}g^1-\dfrac{A_t}{A}\biggr)(g^1+B)
+B_xg^1+2Ag^2_x+B_t=0,
\\[2.5ex]
g^2_t-Ag^2_{xx}-Bg^2_x+\biggl(2g^1_x-\dfrac{A_x}{A}g^1-\dfrac{A_t}{A}\biggr)(g^2-C)-C_xg^1-C_t=0,
\\[2.5ex]
g^3_t-Ag^3_{xx}-Bg^3_x+\biggl(2g^1_x-\dfrac{A_x}{A}g^1-\dfrac{A_t}{A}\biggr)g^3-Cg^3=0,
\end{array}\end{equation}
or an operator $\p_x+\eta(t,x,u)\p_u,$
where the function~$\eta=\eta(t,x,u)$ is a solution of the equation
\begin{equation}\label{EqDet0ForRedOpsOfLPEs}
\eta_t=A(\eta_{xx}+2\eta\eta_{xu}+\eta^2\eta_{uu})+A_x(\eta_x+\eta\eta_u)
+(B\eta)_x+C(\eta-u\eta_u)+C_xu.
\end{equation}
\end{lemma}

\begin{example}
Each equation from class~\eqref{EqGenLPE} with $C=0$ possesses the reduction operator~$\p_x$.
\end{example}

We denote the set of reduction operators of the equation~$\mathcal L$ from class~\eqref{EqGenLPE} by $\mathcal Q(\mathcal L)$, 
omitting the superscript~1. 
The corresponding set factorized with respect to the equivalence of reduction operators is denoted by  $\mathcal Q_{\rm f}(\mathcal L)$. 
Consider the subsets $\mathcal Q_1(\mathcal L)$ and  $\mathcal Q_0(\mathcal L)$ of $\mathcal Q(\mathcal L)$, 
which consist of the operators constrained by the conditions $\tau=1$ and $(\tau,\xi)=(0,1)$, respectively. 
The factor-set $\mathcal Q_{\rm f}(\mathcal L)$ can be identified with $\mathcal Q_1(\mathcal L)\cup\mathcal Q_0(\mathcal L)$. 
This union represents the canonical partition of $\mathcal Q_{\rm f}(\mathcal L)$. 
The systems of form~\eqref{EqDet1ForRedOpsOfLPEs} and equations of form~\eqref{EqDet0ForRedOpsOfLPEs} associated with the equation~$\mathcal L$ 
(and being the determining equations for the operators from $\mathcal Q_1(\mathcal L)$ and  $\mathcal Q_0(\mathcal L)$)
are denoted by ${\rm DE}_1(\mathcal L)$ and ${\rm DE}_0(\mathcal L)$, respectively. 
It is obvious that the rules $\mathcal L\to{\rm DE}_1(\mathcal L)$ and $\mathcal L\to{\rm DE}_0(\mathcal L)$ 
define one-to-one mappings of class~\eqref{EqGenLPE} onto classes~\eqref{EqDet1ForRedOpsOfLPEs} and~\eqref{EqDet0ForRedOpsOfLPEs}.

\begin{note}
The partition of sets of reduction operators according to the condition of \mbox{(non-)}van\-ishing of the coefficient~$\tau$ 
is natural for equations from class~\eqref{EqGenLPE} (as well as the whole class of evolution equations)
and agrees with their transformational properties. 
See Section~\ref{SectionOnLieSymsAndEquivGroupsOfDeterminingEqs} for details.
\end{note}

\begin{note}
For certain reasons, here reduction operators are studied for equations of the non-reduced form~\eqref{EqGenLPE}.
At the same time, it is enough, 
up to the equivalence relation generated by the equivalence group of class~\eqref{EqGenLPE} on the set of pairs 
``(an equation of form~\eqref{EqGenLPE}, its reduction operator)'', to investigate only 
the subclass~\eqref{EqReducedLPE} of equations with $A=1$ and $B=0$. 
The determining equations~\eqref{EqDet1ForRedOpsOfLPEs} and~\eqref{EqDet0ForRedOpsOfLPEs} for 
equations from class~\eqref{EqReducedLPE} have the simpler general form 
\begin{equation}\label{EqReducedDet1ForRedOpsOfLPEs}\arraycolsep=0em
\begin{array}{l}
g^1_t-g^1_{xx}+2g^1_xg^1+2g^2_x=0,
\\[1.5ex]
g^2_t-g^2_{xx}+2g^1_x(g^2+V)+V_xg^1+V_t=0,
\\[1.5ex]
g^3_t-g^3_{xx}+2g^1_xg^3+Vg^3=0
\end{array}\end{equation}
and
\begin{equation}\label{EqReducedDet0ForRedOpsOfLPEs}\arraycolsep=0em
\begin{array}{l}
\eta_t=\eta_{xx}+2\eta\eta_{xu}+\eta^2\eta_{uu}-V(\eta-u\eta_u)-V_xu.
\end{array}\end{equation}
\end{note}

\section{Linearization of determining equations to initial ones}\label{SectionOnLinearizationOfDetEqsForredictionOpsOfLPEs}

There are connections between solution families of an equation~$\mathcal L$ from class~\eqref{EqGenLPE} and its reduction operators. 
This generates connections of the system ${\rm DE}_1(\mathcal L)$ and the equation ${\rm DE}_0(\mathcal L)$ with  
the initial equation~$\mathcal L$ via nonlocal transformations. 

Consider at first reduction operators from $\mathcal Q_1(\mathcal L)$.
Below the indices~$i$ and $j$ run from~1 to~$3$. 
The indices~$p$ and $q$ run from~1 to~$2$. 
The summation convention over repeated indices is used.

\begin{theorem}\label{TheoremUnitedOnSetsOfSolutionsAndReductionOperatorsWithTauNon0OfSystemsOfEqlsp}
Up to the equivalences of operators and solution families, for any equation from the class~\eqref{EqGenLPE} 
there exists a one-to-one correspondence between its reduction operators with nonzero coefficients of $\p_t$
and two-parametric families of its solutions of the form 
\begin{equation}\label{EqLinear2ParamSetOfSolutionsOfLPEs}
u=c_1v^1(t,x)+c_2v^2(t,x)+v^3(t,x),
\end{equation}
where $c_1$ and $c_2$ are constant parameters.
Namely, each operator of such kind corresponds to
the family of solutions which are invariant with respect to this operator. 
The problem of the construction of all two-parametric solution families of an equation~\eqref{EqGenLPE}, 
which are linear in parameters, is completely equivalent 
to the problem of the exhaustive description of its reduction operators with nonzero coefficients of $\p_t$.
\end{theorem}

\begin{corollary}\label{CorollaryOnLinearizationOfDetEqs1ForRedOpsOfLPEs}
Non-linear coupled system~\eqref{EqDet1ForRedOpsOfLPEs} is reduced by the transformation 
\begin{equation}\label{EqTrans1BetweenInitialAndDetEqsForLPEs}\hspace*{-1\arraycolsep}
\begin{array}{l}
g^1=-A\dfrac{v^1v^2_{xx}-v^1_{xx}v^2}{v^1v^2_x-v^1_xv^2}-B,\qquad 
g^2=-A\dfrac{v^1_xv^2_{xx}-v^1_{xx}v^2_x}{v^1v^2_x-v^1_xv^2}+C,\\[3ex] \displaystyle
g^3=\dfrac{A}{v^1v^2_x-v^1_xv^2}
\left|\begin{array}{ccc}
v^1 & v^1_x & v^1_{xx} \\[.2ex]
v^2 & v^2_x & v^2_{xx} \\[.2ex]
v^3 & v^3_x & v^3_{xx}
\end{array}\right|
\end{array}\end{equation}
to the uncoupled system of three copies of equation~\eqref{EqGenLPE} for the functions $v^i=v^i(t,x)$:
\begin{equation}\label{EqUncoupledSystemOf3IdenticalLPEs}
Lv^i=v^i_t-Av^i_{xx}-Bv^i_x-Cv^i=0,
\end{equation}
and the functions $v^1$ and $v^2$ being linearly independent.
\end{corollary}

\begin{note}
Let $W(\varphi^1,\ldots,\varphi^n)$ denote the Wronskian of the functions $\varphi^k=\varphi^k(t,x)$, $k=1,\dots,n$,
with respect to the variable~$x$, i.e. $W(\varphi^1,\ldots,\varphi^n)=\det(\p^l\varphi^k/\p x^l)_{k,l=1}^{\;n}$.
Then transformation~\eqref{EqTrans1BetweenInitialAndDetEqsForLPEs} can be rewritten as 
\[
g^1=-A\dfrac{(W(v^1,v^2))_x}{W(v^1,v^2)}-B,\qquad 
g^2=-A\dfrac{W(v^1_x,v^2_x)}{W(v^1,v^2)}+C,\qquad 
g^3=A\dfrac{W(v^1,v^2,v^3)}{W(v^1,v^2)}.
\]
The solutions $\varphi^k=\varphi^k(t,x)$, $k=1,\dots,n$, of an equation from class~\eqref{EqGenLPE} 
are linearly independent if and only if $W(\varphi^1,\ldots,\varphi^n)\ne0$. See, e.g., Lemma~6 in~\cite{Popovych&Ivanova2004ConsLaws}. 
Therefore, formulas~\eqref{EqTrans1BetweenInitialAndDetEqsForLPEs} are well defined. 
\end{note}

\begin{proof}
Let $\mathcal L$ be an equation from class~\eqref{EqGenLPE} and $Q=\p_t+g^1\p_x+(g^2u+g^3)\p_u\in \mathcal Q_1(\mathcal L)$, i.e., 
the coefficients $g^i=g^i(t,x)$ satisfy the system ${\rm DE}_1(\mathcal L)$. 
An ansatz associated with $Q$ has the form $u=f^1(t,x)\varphi(\omega)+f^0(t,x)$, 
where $f^1=f^1(t,x)$ and $f^0=f^0(t,x)$ are given coefficients, $f^1\ne0$, $\varphi=\varphi(\omega)$ is the new unknown function, 
$\omega=\omega(t,x)$ is the invariant independent variable and $\omega_x\ne0$.
This ansatz reduces $\mathcal L$ to an (in general, inhomogeneous) linear second-order ordinary differential equation in~$\varphi$,
which we denote by $\mathcal L'$. 
The general solution of~$\mathcal L'$ is represented in the form 
$\varphi=c_p\varphi^p(\omega)+\varphi^3(\omega)$, 
where $\varphi^3$ is a particular solution of $\mathcal L'$,  
$\varphi^1$ and $\varphi^2$ are linearly independent solutions of the corresponding homogeneous equation 
and $c_1$ and $c_2$ are arbitrary constants. 
Substituting the general solution of $\mathcal L'$ into the ansatz, 
we obtain the two-parametric family of solutions of~$\mathcal L$, having form~\eqref{EqLinear2ParamSetOfSolutionsOfLPEs} 
with $v^p=f\varphi^p$ and $v^3=f\varphi^3+g$. 
The split in the equations $Lu=0$ and $Q[u]=0$ with respect to the constants~$c_1$ and $c_2$ implies that 
each of the functions $v^i$ is a solution of~$\mathcal L$ and 
\begin{gather*}
(g_1+B)v^p_x-(g^2-C)v^p=-Av^p_{xx}, \\
(g_1+B)v^3_x-(g^2-C)v^3-g^3=-Av^3_{xx}.
\end{gather*}
Since $v^1v^2_x-v^1_xv^2\ne0$, the last system is a well-defined linear system of algebraic equations with respect to $(g^1,g^2,g^3)$, 
whose solution is represented by~\eqref{EqTrans1BetweenInitialAndDetEqsForLPEs}.

Conversely, 
suppose that $\mathcal F$ is a two-parametric family of solutions of~$\mathcal L$, having form~\eqref{EqLinear2ParamSetOfSolutionsOfLPEs}. 
This means that each of the functions $v^i$ is a solution of~$\mathcal L$. 
The functions $v^1$ and $v^2$ are linearly independent since both the parameters $c_1$ and $c_2$ are essential. 
Consider the operator $Q=\p_t+g^1\p_x+(g^2u+g^3)\p_u$, where 
the coefficients $g^i$ are defined by~\eqref{EqTrans1BetweenInitialAndDetEqsForLPEs}. 
$Q[u]=0$ for any $u\in\mathcal F$. 
The ansatz $u=v^1\varphi(\omega)+v^3$, where $\omega=v^2/v^1$, constructed with~$Q$, reduces~$\mathcal L$ 
to the equation $\varphi_{\omega\omega}=0$ since $(v^2/v^1)_x=W(v^1,v^2)/(v^1)^2\ne0$.
Therefore \cite{Zhdanov&Tsyfra&Popovych1999}, $Q\in\mathcal Q_1(\mathcal L)$ and 
the functions $g^i$ have to satisfy the system ${\rm DE}_1(\mathcal L)$. 
\end{proof}

\begin{corollary}\label{CorollaryOnEquivOfRedOps1OfLPEs}
Let $\mathcal L$ be an equation from class~\eqref{EqGenLPE} and 
$G^\infty(\mathcal L)$ denote the trivial Lie invariance group of $\mathcal L$, 
consisting of the linear superposition transformations of the form 
$\tilde t=t$, $\tilde x=x$ and $\tilde u=u+f(t,x)$, 
where the parameter-function $f=f(t,x)$ runs through the solution set of~$\mathcal L$.
Every reduction operator of the equation~$\mathcal L$ with a nonvanishing coefficient of $\p_t$ 
is $G^\infty(\mathcal L)$-equivalent to an operator $\p_t+g^1\p_x+g^2u\p_u$,
where the functions $g^1=g^1(t,x)$ and $g^2=g^2(t,x)$ satisfy the two first equations 
of ${\rm DE}_1(\mathcal L)$. 
\end{corollary}

\begin{proof}
Suppose that a reduction operator~$Q$ of the equation~$\mathcal L$ has a nonvanishing coefficient of $\p_t$. 
In view of Lemma~\ref{LemmaOnReductionOpsWithNonzeroCoeffsOfDtOfLPEs}, 
the operator~$Q$ is equivalent to an operator~$\hat Q$ of the form $\p_t+g^1\p_x+(g^2u+g^3)\p_u,$
where the functions $g^1=g^1(t,x)$, $g^2=g^2(t,x)$ and $g^3=g^3(t,x)$ satisfy the system ${\rm DE}_1(\mathcal L)$.  
It follows from the proof of Theorem~\ref{TheoremUnitedOnSetsOfSolutionsAndReductionOperatorsWithTauNon0OfSystemsOfEqlsp} 
that the coefficient~$g^3$ possesses the representation $g^3=v^3_t+g^1v^3_x-g^2v^3$, 
where $v^3=v^3(t,x)$ is a solution of~$\mathcal L$. 
Then the transformation from $G^\infty(\mathcal L)$ with $f=-v^3$ maps the operator~$\hat Q$ to 
the operator $\tilde Q=\p_t+g^1\p_x+(g^2\tilde u+\tilde g^3)\p_{\tilde u}$, where 
$\tilde g^3=g^3-v^3_t-g^1v^3_x+g^2v^3=0$.
\end{proof}

\begin{note}
The functions $v^i$ satisfying the system~\eqref{EqUncoupledSystemOf3IdenticalLPEs} 
and the additional conditions~\eqref{EqTrans1BetweenInitialAndDetEqsForLPEs} with fixed values of the coefficients $g^j$
are defined up to the transformation 
\begin{equation}\label{EqTransOfViForRedOpsOfLPEs}
\tilde v^p=\mu_{pq}v^q, \quad \tilde v^3=v^3+\mu_{3q}v^q,
\end{equation}
where $\mu_{iq}=\const$, and $\det(\mu_{pq})\ne0$. 
Transformation~\eqref{EqTransOfViForRedOpsOfLPEs} induces the transformation of the constants $c_1$ and $c_2$:
$\tilde c_p=\tilde\mu_{pq}(c_q-\nu_q)$, where $(\tilde\mu_{pq})=(\mu_{p'\!q'\!})^{-1}$. 
It is obvious that the families of solutions~\eqref{EqLinear2ParamSetOfSolutionsOfLPEs} and $u=\tilde c_1\tilde v^1+\tilde c_2\tilde v^2+\tilde v^3$ 
coincides up to re-parametrization and can be identified. 
\end{note}

Consider reduction operators from $\mathcal Q_0(\mathcal L)$.

\begin{theorem}\label{TheoremUnitedOnSetsOfSolutionsAndReductionOperatorsWithTau0OfSystemsOfEqlsp}
Up to the equivalences of operators and solution families, for any equation of form~\eqref{EqGenLPE} 
there exists a one-to-one correspondence between one-parametric families of its solutions 
and reduction operators with zero coefficients of $\p_t$.
Namely, each operator of such kind corresponds to
the family of solutions which are invariant with respect to this operator. 
The problems of the construction of all one-parametric solution families of equation~\eqref{EqGenLPE} 
and the exhaustive description of its reduction operators with zero coefficients of $\p_t$
are completely equivalent.
\end{theorem}

\begin{corollary}\label{CorollaryOnLinearizationOfDetEqs0ForRedOpsOfLPEs}
The nonlinear $(1+2)$-dimensional equation~\eqref{EqDet0ForRedOpsOfLPEs} is reduced by composition of the nonlocal substitution 
$\eta=-\Phi_x/\Phi_u$, where $\Phi$ is a function of $(t,x,u)$, and the hodograph transformation 
\begin{gather}\nonumber
\mbox{the new independent variables:}\qquad\tilde t=t, \quad \tilde x=x, \quad \varkappa=\Phi, 
\\\label{EqTrans0BetweenInitialAndDetEqs}  
\lefteqn{\mbox{the new dependent variable:}}\phantom{\mbox{the new independent variables:}\qquad }\tilde u=u 
\end{gather}
to the initial equation $L\tilde u=0$ in the function $\tilde u=\tilde u(\tilde t,\tilde x,\varkappa)$ 
with $\varkappa$ playing the role of a parameter.
\end{corollary}

\begin{proof}
Let $\mathcal L$ be an equation from class~\eqref{EqGenLPE} and $Q=\p_x+\eta\p_u\in \mathcal Q_0(\mathcal L)$, i.e., 
the coefficient $\eta=\eta(t,x,u)$ satisfies the equation ${\rm DE}_0(\mathcal L)$. 
An ansatz associated with $Q$ has the form $u=f(t,x,\varphi(\omega))$, 
where $f=f^1(t,x,\varphi)$ is a given function, $f_\varphi\ne0$, 
$\varphi=\varphi(\omega)$ is the new unknown function and $\omega=t$ is the invariant independent variable.
This ansatz reduces $\mathcal L$ to a first-order ordinary differential equation in~$\varphi$, which we denote by $\mathcal L'$. 
The general solution of~$\mathcal L'$ is represented in the form 
$\varphi=\varphi(\omega,\varkappa)$, where $\varphi_\varkappa\ne0$ and $\varkappa$ is an arbitrary constant. 
The substitution of the general solution of $\mathcal L'$ into the ansatz results in  
the one-parametric family~$\mathcal F$ of solutions $u=\tilde f(t,x,\varkappa)$ of~$\mathcal L$ 
with $\tilde f=f(t,x,\varphi(t,\varkappa))$. 
Expressing the parameter~$\varkappa$ from the equality  $u=\tilde f(t,x,\varkappa)$, 
we obtain that $\varkappa=\Phi(t,x,u)$, where $\Phi_u\ne0$. 
Then $\eta=u_x=-\Phi_x/\Phi_u$ for any $u\in\mathcal F$, i.e., for any admissible values of $(t,x,\varkappa)$. 
This implies that $\eta=-\Phi_x/\Phi_u$ for any admissible values of $(t,x,u)$.

Conversely, 
suppose that $\mathcal F=\{u=f(t,x,\varkappa)\}$ is a one-parametric family of solutions of~$\mathcal L$. 
The derivative $f_\varkappa$ is nonzero since the parameter $\varkappa$ is essential. 
We express $\varkappa$ from the equality $u=f(t,x,\varkappa)$: $\varkappa=\Phi(t,x,u)$ for some function $\Phi=\Phi(t,x,u)$ 
with $\Phi_u\ne0$. 
Consider the operator $Q=\p_x+\eta\p_u$, where the coefficient $\eta=\eta(t,x,u)$ is defined by the formula $\eta=-\Phi_x/\Phi_u$. 
$Q[u]=0$ for any $u\in\mathcal F$. 
The ansatz $u=f(t,x,\varphi(\omega))$, where $\omega=t$, associated with~$Q$, 
reduces~$\mathcal L$ to the equation $\varphi_\omega=0$.
Therefore \cite{Zhdanov&Tsyfra&Popovych1999}, $Q\in\mathcal Q_0(\mathcal L)$ and hence
the function $\eta$ satisfies the equation ${\rm DE}_0(\mathcal L)$. 
\end{proof}

\begin{note}\label{NoteOnEquiv1ParamFamiliesOfSolutions}
One-parametric families of solutions $u=f(t,x,\varkappa)$ and $u=\tilde f(t,x,\tilde\varkappa)$ of~$\mathcal L$ 
are assumed equivalent if they consist of the same solutions and differ only by parameterizations, i.e., 
if there exists a function $\zeta=\zeta(\varkappa)$ such that 
$\zeta_\varkappa\ne0$ and $\tilde f(t,x,\zeta(\varkappa))=f(t,x,\varkappa)$.  
Equivalent one-parametric families of solutions are associated with the same operator 
from $\mathcal Q_0(\mathcal L)$ and have to be identified. 
\end{note}

\begin{note}
The supposed triviality of the above ansatzes and reduced equations is connected with 
usage of the special representations for the solutions of the determining equations.  
Under this approach difficulties in the construction of ansatzes and the integration of reduced equations are replaced 
by difficulties in obtaining the representations for coefficients of reduction operators.
\end{note}

\section{Admissible transformations, the equivalence groups\\ and Lie symmetries of determining equations}
\label{SectionOnLieSymsAndEquivGroupsOfDeterminingEqs}

The ``no-go'' results of the previous section can be extended with investigation of point transformations, Lie symmetries and Lie reductions 
of determining equations~\eqref{EqDet1ForRedOpsOfLPEs} and~\eqref{EqDet0ForRedOpsOfLPEs}. 
Thus, the maximal Lie invariance algebras of~\eqref{EqDet1ForRedOpsOfLPEs} 
and~\eqref{EqDet0ForRedOpsOfLPEs} are isomorphic 
to the maximal Lie invariance algebras of equation~\eqref{EqGenLPE} in a canonical way. 
(Before this result was known only for the linear heat equation~\cite{Fushchych&Shtelen&Serov&Popovych1992}.) 
Moreover, the similar statements are true for the complete point symmetry groups including discrete symmetry transformations as well as  
the equivalence groups and sets of admissible transformations of classes of the above equations. 

All these statements are justified by Lemmas~\ref{LemmaOnInducedMappingOfRedictionOps} and~\ref{LemmaOnAdmTransOfInhomLPEs}. 
Indeed, each point transformation $\mathcal T$ between equations $\mathcal L$ and $\tilde{\mathcal L}$ from class~\eqref{EqGenLPE} 
has the form~\eqref{EqGenFormOfTransOfInhomLPEs} and induces the one-to-one mappings 
$\mathcal T_*\colon\mathcal Q(\mathcal L)\to\mathcal Q(\tilde{\mathcal L})$ and 
$\mathcal T_{\rm f}\colon\mathcal Q_{\rm f}(\mathcal L)\to\mathcal Q_{\rm f}(\tilde{\mathcal L})$. 
Due to the conditions $\mathcal T^t_x=0$ and $\mathcal T^t_u=0$, the transformation $\mathcal T_*$ preserves the constraint $\tau=0$ 
(resp. $\tau\ne0$) for coefficients of reduction operators.
Therefore, the transformation~$\mathcal T_{\rm f}$ is split into the one-to-one mappings 
$\mathcal T_{{\rm f},1}\colon\mathcal Q_1(\mathcal L)\to\mathcal Q_1(\tilde{\mathcal L})$ and 
$\mathcal T_{{\rm f},0}\colon\mathcal Q_0(\mathcal L)\to\mathcal Q_0(\tilde{\mathcal L})$ 
according to the canonical partitions of $\mathcal Q_{\rm f}(\mathcal L)$ and $\mathcal Q_{\rm f}(\tilde{\mathcal L})$.
This implies that there exist the transformations~$\mathcal T_1$ and~$\mathcal T_0$ in the spaces of 
the variables $(t,x,g^1,g^2,g^3)$ and $(t,x,u,\eta)$, 
which are induced by the transformation~$\mathcal T$ in a canonical way. 
It is evident that 
\[
\mathcal T_1\bigl({\rm DE}_1(\mathcal L)\bigr)={\rm DE}_1(\tilde{\mathcal L}),\quad
\mathcal T_0\bigl({\rm DE}_0(\mathcal L)\bigr)={\rm DE}_0(\tilde{\mathcal L}).
\]

The procedure of deriving the explicit formulas for $\mathcal T_1$ is the following: 
Acting on the operator $\p_t+g^1\p_x+(g^2u+g^3)\p_u$ by $\mathcal T_*$ and then 
normalizing the coefficient of~$\p_{\tilde t}$ to~1, we obtain the operator  
$\p_{\tilde t}+\tilde g^1\p_{\tilde x}+(\tilde g^2\tilde u+\tilde g^3)\p_{\tilde u}$, 
where the new coefficients $\tilde g^i=\tilde g^i(\tilde t,\tilde x)$, $i=1,2,3$, are calculated by the formulas 
\begin{gather}\nonumber
\tilde g^1=\frac{X_x}{T_t}g^1+\frac{X_t}{T_t}, \\\label{EqTransOfCoeffsOfRedOps1ForLPEs}
\tilde g^2=\frac{1}{T_t}g^2+\frac{U^1_x}{T_tU^1}g^1+\frac{U^1_t}{T_tU^1}, \\\nonumber
\tilde g^3=\frac{U^1}{T_t}g^3-\frac{U^0}{T_t}g^2+\frac{U^0_xU^1-U^0U^1_x}{T_tU^1}g^1+\frac{U^0_tU^1-U^0U^1_t}{T_tU^1}. 
\end{gather}
Formulas~\eqref{EqTransOfCoeffsOfRedOps1ForLPEs} describe the action of~$\mathcal T_1$ on the dependent variables $(g^1,g^2,g^3)$. 
The independent variables $t$ and $x$ and the arbitrary elements $A$, $B$ and $C$ are transformed 
by the same formulas~\eqref{EqGenFormOfTransOfInhomLPEs} and~\eqref{EqTransOfCoeffsOfLPE} as ones of the transformation~$\mathcal T$.
The transformation of~$u$ is neglected.

If the transformation~$\mathcal T$ belongs to the equivalence group~$G^\sim$ of class~\eqref{EqGenLPE} 
then it is defined for all values of arbitrary elements. 
Therefore, the same statement is true for~$\mathcal T_1$, i.e., 
$\mathcal T_1$ belongs to the equivalence group~$G^\sim_1$ of class~\eqref{EqDet1ForRedOpsOfLPEs}. 
In other words, the equivalence group of the initial class induces a subgroup of the equivalence group 
of the class of determining equations for the case $\tau=1$.

Suppose that the transformation~$\mathcal T$ is parameterized by the parameter $\varepsilon$ 
and this family of transformations form a one-parametric Lie symmetry group of the equation~$\mathcal L$, generated by an operator 
$Q=\tau\p_t+\xi\p_x+(\zeta^1u+\zeta^0)\p_u$. 
We differentiate formulas~\eqref{EqTransOfCoeffsOfRedOps1ForLPEs} with respect to $\varepsilon$ and then 
put $\varepsilon=0$, taking into account the conditions 
\[\arraycolsep=0ex
\begin{array}{llll}
\tau=\tau(t)=T_\varepsilon\bigl|_{\varepsilon=0},\quad &T\bigl|_{\varepsilon=0}=t,\quad
&\xi=\xi(t,x)=X_\varepsilon\bigl|_{\varepsilon=0},\quad &X\bigl|_{\varepsilon=0}=x,
\\[1.5ex]
\zeta^1=\zeta^1(t,x)=U^1_\varepsilon\bigl|_{\varepsilon=0},\quad &U^1\bigl|_{\varepsilon=0}=1,\quad
&\zeta^0=\zeta^0(t,x)=U^0_\varepsilon\bigl|_{\varepsilon=0},\quad &U^0\bigl|_{\varepsilon=0}=0.
\end{array}\]
As a result, we obtain the expressions for the coefficients $\theta^i$ 
of the Lie symmetry operator $Q_1=\tau\p_t+\xi\p_x+\theta^i\p_{g^i}$ of the system ${\rm DE}_1(\mathcal L)$, 
associated with the operator~$Q$:
\begin{gather}\nonumber
\theta^1=(\xi_x-\tau_t)g^1+\xi_t, \\\label{EqCoeffsOfSymOpersOfDetSystemOfRedOps1ForLPEs}
\theta^2=-\tau_tg^2+\eta^1_xg^1+\eta^1_t, \\\nonumber
\theta^3=(\eta^1-\tau_t)g^3-\eta^0g^2+\eta^0_xg^1+\eta^0_t. 
\end{gather}

The explicit formulas for $\mathcal T_0$ are derived in the analogous way.
The action of $\mathcal T_*$ on the operator $\p_x+\eta\p_u$ and 
the normalization of the coefficient of~$\p_{\tilde x}$ to~1 result in the operator $\p_{\tilde x}+\tilde\eta\p_{\tilde u}$, where 
\begin{equation}\label{EqTransOfCoeffsOfRedOps0ForLPEs}
\tilde\eta=\frac{U^1}{X_x}\eta+\frac{U^1_x}{X_x}u+\frac{U^0_x}{X_x}.
\end{equation}
Formula~\eqref{EqTransOfCoeffsOfRedOps0ForLPEs} represents the expression for the dependent variable $\eta$ transformed by~$\mathcal T_0$.
The transformations of independent variables $t$, $x$ and $u$ and the arbitrary elements $A$, $B$ and $C$ 
are given by the formulas~\eqref{EqGenFormOfTransOfInhomLPEs} and~\eqref{EqTransOfCoeffsOfLPE}. 
The unique difference from the transformation~$\mathcal T$ is that the variable $u$ is assumed independent.
This implies that each transformation from the equivalence group~$G^\sim$ of class~\eqref{EqGenLPE} induces 
a transformation from the equivalence group~$G^\sim_0$ of class~\eqref{EqDet0ForRedOpsOfLPEs}.

Under the infinitesimal approach, each Lie invariance operator $Q=\tau\p_t+\xi\p_x+(\zeta^1u+\zeta^0)\p_u$ 
of~$\mathcal L$ is prolonged to the Lie invariance operator $Q_0=Q+\theta\p_\eta$ of ${\rm DE}_0(\mathcal L)$, 
where the coefficient~$\theta$ is determined by the formula
\begin{equation}\label{EqCoeffsOfSymOpersOfDetSystemOfRedOps0ForLPEs}
\theta=(\zeta^1-\xi_x)\eta+\zeta^1_xu+\zeta^0_x. 
\end{equation}

The problem is to prove that the induced objects 
(resp. admissible transformations, point equivalences, point symmetries and Lie invariance operators)
exhaust all possible objects of the corresponding kinds for determining equations. 

\begin{lemma}\label{LemmaOnAdmTransOfDetSystemOfRedOps1ForLPEs}
If a point transformation connects two systems ${\rm DE}_1(\mathcal L)$ and ${\rm DE}_1(\tilde{\mathcal L})$ 
from class~\eqref{EqDet1ForRedOpsOfLPEs} then it has the form 
\begin{equation}\label{EqGenFormOfdmTransOfDetSystemOfRedOps1ForLPEs-1}
\tilde t=T(t), \quad \tilde x=X(t,x), \quad \tilde g^i=G^{ii'\!}(t,x)g^{i'\!}+G^{i0}(t,x),
\end{equation}
where $T$, $X$, $G^{33}$ and $G^{32}$ are smooth functions of their arguments such that $T_tX_xG^{33}\ne0$ 
and additionally $G^{32}/G^{33}$ is a solution of the associated equation~$\mathcal L$; $i,i'=1,2,3$. 
The other parameter-functions in~\eqref{EqGenFormOfdmTransOfDetSystemOfRedOps1ForLPEs-1} are explicitly defined: 
\begin{equation}\label{EqGenFormOfdmTransOfDetSystemOfRedOps1ForLPEs-2}\arraycolsep=0ex
\begin{array}{llll}
G^{10}=\dfrac{X_t}{T_t},\qquad& G^{11}=\dfrac{X_x}{T_t},\qquad& G^{12}=0,\qquad& G^{13}=0,
\\[2.5ex]
G^{20}=\dfrac{(T_tG^{33})_t}{T_t{}^2G^{33}},\qquad& G^{21}=\dfrac{G^{33}_x}{T_tG^{33}},\qquad& G^{22}=\dfrac1{T_t},\qquad& G^{23}=0,
\\[2.5ex]
G^{30}=\dfrac{(T_tG^{33})_t}{T_t{}^2G^{33}},\qquad& G^{31}=\dfrac{G^{33}_x}{G^{33}}G^{32}-G^{32}_x.& \qquad& 
\end{array}
\end{equation}
The arbitrary elements are transformed by the formulas
\begin{gather}\label{EqGenFormOfdmTransOfDetSystemOfRedOps1ForLPEs-3}
\tilde A=\frac{X_x^2}{T_t}A,\quad 
\tilde B=\frac{X_x}{T_t}\left(B-2\frac{G^{33}_x}{G^{33}}A\right)-\frac{X_t-AX_{xx}}{T_t},\quad 
\tilde C=-G^{33}L\frac1{T_tG^{33}}.
\end{gather}
Here, $L=\p_t-A\p_{xx}-B\p_x-C$ is the second-order linear differential operator associated with 
the equation~$\mathcal L$.
\end{lemma}

\begin{proof}
The systems ${\rm DE}_1(\mathcal L)$ and ${\rm DE}_1(\tilde{\mathcal L})$ consists of 
second-order evolution equations which are linear in the derivatives, and coefficients of 
second derivatives form the nonsingular matrices $\diag(A,A,A)$ and $\diag(\tilde A,\tilde A,\tilde A)$, respectively. 
In view of Corollary~13 of~\cite{Popovych&Kunzinger&Ivanova2007} each transformation between such systems 
necessarily has form~\eqref{EqGenFormOfdmTransOfDetSystemOfRedOps1ForLPEs-1}.
We apply the direct method with taking into account the conditions~\eqref{EqGenFormOfdmTransOfDetSystemOfRedOps1ForLPEs-1} 
and find more conditions which can be split by $g^i$ and $g^i_x$. 
The system of determining equations on parameters of the transformation, obtained after the split, 
implies equations~\eqref{EqGenFormOfdmTransOfDetSystemOfRedOps1ForLPEs-2} and 
expressions~\eqref{EqTransOfCoeffsOfLPE} for transformations of the arbitrary elements. 
\end{proof}

\begin{theorem}\label{TheoremOnAdmTransOfDetSystemOfRedOps1ForLPEs}  
There exists a canonical one-to-one correspondence between the sets of admissible transformations 
of classes~\eqref{EqGenLPE} and~\eqref{EqDet1ForRedOpsOfLPEs}. 
Namely, each point transformation between equations $\mathcal L$ and $\tilde{\mathcal L}$ from class~\eqref{EqGenLPE} 
induces a point transformation between the associated systems ${\rm DE}_1(\mathcal L)$ and ${\rm DE}_1(\tilde{\mathcal L})$ 
according to formulas~\eqref{EqTransOfCoeffsOfRedOps1ForLPEs}. 
In both the transformations the independent variables are transformed in the same way. 
The induced transformations exhaust the sets of admissible transformation in class~\eqref{EqDet1ForRedOpsOfLPEs}.
\end{theorem} 

\begin{proof}\looseness=-1
It  only remains to prove that every admissible transformation in class~\eqref{EqDet1ForRedOpsOfLPEs} is induced by 
an admissible transformation in class~\eqref{EqGenLPE} in the above way. 
We fix two point-equivalent systems from class~\eqref{EqDet1ForRedOpsOfLPEs}. 
They necessarily are systems of determining equations for reduction operators with the unit coefficients of $\p_t$ 
for some equations $\mathcal L$ and $\tilde{\mathcal L}$ from class~\eqref{EqGenLPE}. 
Therefore, these systems can be denoted by ${\rm DE}_1(\mathcal L)$ and ${\rm DE}_1(\tilde{\mathcal L})$, respectively.
Consider a point transformation~$\breve{\mathcal T}$ mapping the system ${\rm DE}_1(\mathcal L)$ to the system ${\rm DE}_1(\tilde{\mathcal L})$. 
In view of Lemma~\ref{LemmaOnAdmTransOfDetSystemOfRedOps1ForLPEs}, the transformation~$\breve{\mathcal T}$ 
has form~\eqref{EqGenFormOfdmTransOfDetSystemOfRedOps1ForLPEs-1}, where $G^{32}/G^{33}$ is a solution of~$\mathcal L$ 
and the other parameter-functions $G^{ii'}$ and $G^{i0}$ are explicitly expressed by~\eqref{EqGenFormOfdmTransOfDetSystemOfRedOps1ForLPEs-2}.
Formulas~\eqref{EqGenFormOfdmTransOfDetSystemOfRedOps1ForLPEs-3} describe connections between the arbitrary elements 
of ${\rm DE}_1(\mathcal L)$ and ${\rm DE}_1(\tilde{\mathcal L})$. 
We  associate the transformation~$\breve{\mathcal T}$ with the transformation~$\mathcal T$ in the space of the variables $(t,x,u)$, 
having the form~\eqref{EqGenFormOfTransOfInhomLPEs}, where $U^1=T_tG^{33}$ and $U^0=T_tG^{32}$.
By the construction, $U^1/U^0$ is a solution of~$\mathcal L$.
Since the pairs $({\rm DE}_1(\mathcal L),{\rm DE}_1(\tilde{\mathcal L}))$ and $(\mathcal L,\tilde{\mathcal L})$ have the same tuples of 
arbitrary elements, Lemma~\ref{LemmaOnAdmTransOfInhomLPEs} and formulas~\eqref{EqGenFormOfdmTransOfDetSystemOfRedOps1ForLPEs-3} imply 
that $\mathcal T$ is a point transformation from $\mathcal L$ to $\tilde{\mathcal L}$. 
The comparison of~\eqref{EqGenFormOfdmTransOfDetSystemOfRedOps1ForLPEs-2} with~\eqref{EqTransOfCoeffsOfRedOps1ForLPEs} 
allows us to conclude that $\breve{\mathcal T}$ is induced by~$\mathcal T$, i.e., $\breve{\mathcal T}=\mathcal T_1$.  
\end{proof}

\begin{note}
It follows from the proof of Theorem~\ref{TheoremOnAdmTransOfDetSystemOfRedOps1ForLPEs} that 
``if \dots\ then \dots'' in Lemma~\ref{LemmaOnAdmTransOfDetSystemOfRedOps1ForLPEs} can be replaced by ``\dots\ if and only if \dots'', i.e., 
the presented conditions are necessary and sufficient.   
\end{note}

\begin{corollary}\label{CorollaryOnEquivGroupOfDetSystemOfRedOps1ForLPEs}
The equivalence group $G^\sim_1$ of class~\eqref{EqDet1ForRedOpsOfLPEs} is isomorphic to the equivalence group $G^\sim$ of class~\eqref{EqGenLPE}. 
The canonical isomorphism is established by formulas~\eqref{EqTransOfCoeffsOfRedOps1ForLPEs}, where $U^0=0$. 
\end{corollary}

\begin{corollary}\label{CorollaryOnLieSymAlgOfDetSystemOfRedOps1ForLPEs}
For each equation~$\mathcal L$ from class~\eqref{EqGenLPE}, the maximal point symmetry groups 
(resp. the maximal Lie invariance algebras) of the equation~$\mathcal L$ and the system ${\rm DE}_1(\mathcal L)$ are isomorphic. 
A Lie symmetry operator $Q=\tau\p_t+\xi\p_x+(\zeta^1u+\zeta^0)\p_u$ of~$\mathcal L$ induces the 
the Lie symmetry operator $Q_1=\tau\p_t+\xi\p_x+\theta^i\p_{g^i}$ of the system ${\rm DE}_1(\mathcal L)$, 
where the coefficients $\theta^i$, $i=1,2,3$, are defined by formulas~\eqref{EqCoeffsOfSymOpersOfDetSystemOfRedOps1ForLPEs}. 
\end{corollary}

Corollaries~\ref{CorollaryOnEquivGroupOfDetSystemOfRedOps1ForLPEs} and \ref{CorollaryOnLieSymAlgOfDetSystemOfRedOps1ForLPEs} 
along with Theorem~\ref{TheoremOnGroupClassificationOfLPEs} give the group classification of class~\eqref{EqDet1ForRedOpsOfLPEs}.
  
\begin{corollary}\label{CorollaryOnGroupClassificationOfDetSystemOfRedOps1ForLPEs}
The kernel Lie algebra of class~\eqref{EqDet1ForRedOpsOfLPEs} is $\langle I_1\rangle$, where $I_1=g^3\p_{g^3}$.
Any system from class~\eqref{EqDet1ForRedOpsOfLPEs} is invariant with respect to the operators of the form $Z_1(f)=(f_t+f_xg^1-fg^2)\p_{g^3}$,
where the parameter-function $f=f(t,x)$ runs through the solution set of the associated equation 
$f_t=Af_{xx}+Bf_x+Cf$.
All possible $G^\sim_1$-inequivalent cases of extension of the maximal 
Lie invariance algebra are exhausted by the following systems of the reduced form~\eqref{EqReducedDet1ForRedOpsOfLPEs}  
(the values of~$V$ are given together with the corresponding maximal Lie invariance algebras):

\vspace{1.5ex}

$\makebox[6mm][l]{\rm 1.}
V=V(x) \colon\quad \langle\p_t,\,I_1,\,Z_1(f)\rangle;$

\vspace{1ex}

$\makebox[6mm][l]{\rm 2.}
V=\mu x^{-2},\ \mu\ne0 \colon\quad \langle\p_t,\, D_1,\, \Pi_1,\, I_1,\,Z_1(f)\rangle;$

\vspace{1ex}

$\makebox[6mm][l]{\rm 3.}
V=0 \colon\quad \langle\p_t,\, \p_x,\, G_1,\, D_1,\, \Pi_1,\, I_1,\,Z_1(f)\rangle$.

\vspace{1.5ex}

\noindent
Here 
\begin{gather*}
D_1=2t\p_t+x\p_x-g^1\p_{g^1}-2g^2\p_{g^2}, 
\\[.5ex]
\Pi_1=4t^2\p_t+4tx\p_x+4(x-tg^1)\p_{g^1}-(8tg^2+2xg^1+2)\p_{g^2}-(x^2+10t)g^3\p_{g^3}, 
\\[.5ex]
G_1=2t\p_x+2\p_{g^1}-g^1\p_{g^2}-xg^3\p_{g^3}.
\end{gather*}
\end{corollary}

\begin{note}
It is obvious that 
Corollaries~\ref{CorollaryOnEquivGroupOfDetSystemOfRedOps1ForLPEs}, \ref{CorollaryOnLieSymAlgOfDetSystemOfRedOps1ForLPEs} 
and~\ref{CorollaryOnGroupClassificationOfDetSystemOfRedOps1ForLPEs} can be reformulated 
for subclass~\eqref{EqReducedLPE} of the initial equations in the reduced form 
and subclass~\eqref{EqReducedDet1ForRedOpsOfLPEs} of the corresponding determining equations of the first kind (the case $\tau\ne0$).
\end{note}

A specific question for class~\eqref{EqDet1ForRedOpsOfLPEs} is 
what transformations of the functions $(v^1,v^2,v^3)$ defined in Corollary~\ref{CorollaryOnLinearizationOfDetEqs1ForRedOpsOfLPEs} 
are induced by admissible transformations in class~\eqref{EqDet1ForRedOpsOfLPEs}.
It is clear that each induced transformation is admissible in class~\eqref{EqUncoupledSystemOf3IdenticalLPEs}. 
Let $\mathcal L$ and $\tilde{\mathcal L}$ be equations from class~\eqref{EqGenLPE}. 
Denote the corresponding systems of form~\eqref{EqUncoupledSystemOf3IdenticalLPEs} by $3\mathcal L$ and $3\tilde{\mathcal L}$ 
and the corresponding sets of formulas~\eqref{EqTrans1BetweenInitialAndDetEqsForLPEs} by $\mathcal G$ and $\tilde{\mathcal G}$, 
respectively.
It is proved analogously to Lemma~\ref{LemmaOnAdmTransOfDetSystemOfRedOps1ForLPEs} that
any point transformation connecting the systems $3\mathcal L$ and $3\tilde{\mathcal L}$ has the form 
\[
\tilde t=T(t), \quad \tilde x=X(t,x), \quad \tilde v^i=U^1(t,x)\mu_{ij}v^j+U^{i0}(t,x),
\]
where $\mu_{ij}=\const$, $\det(\mu_{ij})\ne0$, $i,j=1,2,3$; 
$T$, $X$, $U^1$ and $U^{i0}$ are arbitrary smooth functions of their arguments such that $T_tX_xU^1\ne0$ 
and additionally $U^{i0}/U^1$ are solutions of the equation~$\mathcal L$. 
The arbitrary elements are transformed by the formulas~\eqref{EqTransOfCoeffsOfLPE}, 
where $L=\p_t-A\p_{xx}-B\p_x-C$ is the second-order linear differential operator associated with the equation~$\mathcal L$.
The agreement of transformations between $3\mathcal L$ and $3\tilde{\mathcal L}$ with transformations between 
${\rm DE}_1(\mathcal L)$ and ${\rm DE}_1(\tilde{\mathcal L})$ via formulas~\eqref{EqTrans1BetweenInitialAndDetEqsForLPEs}
implies the additional conditions 
\[
\mu_{13}=\mu_{23}=0, \quad U^{10}=U^{20}=0, \quad \mu_{33}=1, \quad U^1=T_tG^{33}, \quad U^{30}=T_tG^{30},
\] 
for the admissible transformations between the systems 
$3\mathcal L\cap\mathcal G\cap{\rm DE}_1(\mathcal L)$ and $3\tilde{\mathcal L}\cap\tilde{\mathcal G}\cap{\rm DE}_1(\tilde{\mathcal L})$.
To derive these conditions, we express all the tilde variables in~$\tilde{\mathcal G}$ via the non-tilde ones, 
then substitute the expressions for $g^i$ given by~$\mathcal G$ into~$\tilde{\mathcal G}$ and split with respect to 
$v^j$ and their derivatives.  
Combining the obtained result with Theorem~\ref{TheoremOnAdmTransOfDetSystemOfRedOps1ForLPEs} 
and omitting the systems ${\rm DE}_1(\mathcal L)$ and ${\rm DE}_1(\tilde{\mathcal L})$ as differential consequences of 
the systems $3\mathcal L\cap\mathcal G$ and $3\tilde{\mathcal L}\cap\tilde{\mathcal G}$, respectively, we get 
that the point transformation~$\mathcal T$ of form~\eqref{EqGenFormOfTransOfInhomLPEs} 
between the equations $\mathcal L$ and $\tilde{\mathcal L}$ induces the point transformation 
\[
\tilde t=T(t), \quad \tilde x=X(t,x), \quad \tilde v^p=U^1(t,x)\mu_{pq}v^q, \quad \tilde v^3=U^1(t,x)\mu_{3q}v^q+U^0(t,x),
\]
where $\det(\mu_{pq})\ne0$, $p,q=1,2$, between the system $3\mathcal L\cap\mathcal G$ and $3\tilde{\mathcal L}\cap\tilde{\mathcal G}$.
The appearance of the additional constants~$\mu_{iq}$ in the induced transformation is explained by 
the uncertainty~\eqref{EqTransOfViForRedOpsOfLPEs} under determining the function~$v^i$.
The consideration of a one-parametric Lie symmetry group of the equation~$\mathcal L$ instead of a single transformation 
between the (possibly different) equations $\mathcal L$ and $\tilde{\mathcal L}$ results in 
a formula for the extension of Lie symmetry operators of~$\mathcal L$ to Lie symmetry operators of~$3\mathcal L$. 
Namely, the following statement is true.

\begin{lemma}\label{LemmaOnLieSymOpsOfLinearizedDetSystemOfRedOps1ForLPEs}
Each Lie symmetry operator $Q=\tau\p_t+\xi\p_x+(\zeta^1u+\zeta^0)\p_u$ of the equation~$\mathcal L$ generates  
the family 
\[
\{\tau\p_t+\xi\p_x+\zeta^1v^i\p_{v^i}+\zeta^0\p_{v^3}+\lambda_{iq}v^q\p_{v^i}\mid \lambda_{iq}=\const\}
\]
of Lie symmetry operators of the associated system~$3\mathcal L$ with the additional conditions~$\mathcal G$. 
Here $i,j=1,2,3$, $q=1,2$.
The functions $g^j$ satisfy the system~${\rm DE}_1(\mathcal L)$ being the compatibility condition of $3\mathcal L\cap\mathcal G$. 
\end{lemma}

The chain of similar statements is also obtained for class~\eqref{EqDet0ForRedOpsOfLPEs}.

\begin{lemma}\label{LemmaOnAdmTransOfDetSystemOfRedOps0ForLPEs}
If a point transformation in the space of the variables $(t,x,u,\eta)$ connects 
two equations ${\rm DE}_0(\mathcal L)$ and ${\rm DE}_0(\tilde{\mathcal L})$ 
from class~\eqref{EqDet0ForRedOpsOfLPEs} then it has the form given 
by formulas~\eqref{EqGenFormOfTransOfInhomLPEs} and~\eqref{EqTransOfCoeffsOfRedOps0ForLPEs}, 
where $T$, $X$, $U^1$ and $U^0$ are arbitrary smooth functions of their arguments such that $T_tX_xU^1\ne0$ 
and additionally $U^0/U^1$ is a solution of the equation~$\mathcal L$. 
The arbitrary elements are transformed by the formulas~\eqref{EqTransOfCoeffsOfLPE}, 
where $L=\p_t-A\p_{xx}-B\p_x-C$ is the second-order linear differential operator associated with the equation~$\mathcal L$.
\end{lemma}

\begin{proof}
The matrices formed by the coefficients of the second derivations in the equations ${\rm DE}_0(\mathcal L)$ and ${\rm DE}_0(\tilde{\mathcal L})$ 
are singular. That is why we cannot use the results of~\cite{Prokhorova2005} on admissible transformations in classes of parabolic equations 
having positively defined matrices of the coefficients of the second derivations. 
All determining equations have to be obtained independently. 

We use the direct method. 
Consider a point transformation $\mathcal T$ from the equation ${\rm DE}_0(\mathcal L)$ to the equation ${\rm DE}_0(\tilde{\mathcal L})$ 
of the general form $[\tilde t,\tilde x,\tilde u,\tilde\eta]=[T,X,U,H](t,x,u,\eta)$ with the nonvanishing Jacobian. 
Sometimes we will also assume that the old variables $(t,x,u,\eta)$ are functions of the new variables $(\tilde t,\tilde x,\tilde u,\tilde\eta)$ 
and do a simultaneous split with respect to both the old and new variables. This trick is correct under certain conditions. 
We introduce the notations $Q:=D_x+\eta D_u$, $\tilde Q:=D_{\tilde x}+\tilde\eta D_{\tilde u}$ and $F:=\tilde Q\tilde\eta$. 
In the old variables, the function $F$ is expressed via $t$, $x$, $u$, $\eta$, $\eta_t$, $\eta_x$ and $\eta_u$, 
and moreover $(F_{\eta_t},F_{\eta_x},F_{\eta_u})\ne(0,0,0)$. 
(Indeed, the condition $F_{\eta_t}=F_{\eta_x}=F_{\eta_u}=0$ means that the function $F$ depends only on 
$(t,x,t,\eta)$ in the old variables and, therefore, is a function of only $(\tilde t,\tilde x,\tilde u,\tilde\eta)$ in the new variables. 
Then we could split the equation $F=\tilde\eta_{\tilde x}+\tilde\eta\tilde\eta_{\tilde u}$ defining~$F$
with respect to derivatives of~$\tilde\eta$ and obtain the contradiction $0=1$.)

The equation ${\rm DE}_0(\tilde{\mathcal L})$ can be written in the form $\tilde QF=\dots$, 
where the right-hand side contains derivatives only up to order~1.
We return to the old variables in ${\rm DE}_0(\tilde{\mathcal L})$ and confine it to the manifold 
of the equation ${\rm DE}_0(\mathcal L)$, expressing the derivative $\eta_{xx}$ from ${\rm DE}_0(\mathcal L)$
and substituting the found expression into ${\rm DE}_0(\tilde{\mathcal L})$. 
Then we split the obtained equation ${\rm DE}_0'$ step by step with respect to different subsets of the other derivatives of $\eta$ 
(or $\tilde\eta$ alternatively). 
To optimize this procedure, we start from the subsets of derivatives giving the simplest determining equations 
and take into account found equations for the further split. 
Note that the expression $\tilde QF$ has the representation $\tilde QF=(\tilde Qt)D_tF+(\tilde Qx)D_xF+(\tilde Qu)D_uF$. 

After collecting the coefficients of $\eta_{tt}$, $\eta_{tx}$ and $\eta_{tu}$ in ${\rm DE}_0'$, 
we derive the system
\[
(\tilde Qt)F_{\eta_t}=0, \quad
(\tilde Qt)F_{\eta_x}+(\tilde Qx)F_{\eta_t}=0, \quad
(\tilde Qt)F_{\eta_u}+(\tilde Qu)F_{\eta_t}=0 
\]
which implies the equation $\tilde Qt=0$ since $(F_{\eta_t},F_{\eta_x},F_{\eta_u})\ne(0,0,0)$. 
We expand the expression $\tilde Qt$, assuming $t$ a function of $(\tilde t,\tilde x,\tilde u,\tilde\eta)$: 
$\tilde Qt=t_{\tilde x}+t_{\tilde\eta}\tilde\eta_{\tilde x}+\tilde\eta(t_{\tilde u}+t_{\tilde\eta}\tilde\eta_{\tilde u})$.  
The split of the equation $\tilde Qt=0$ with respect to the new jet variables $\tilde\eta_{\tilde x}$ and $\tilde\eta_{\tilde u}$ 
results in the equations $t_{\tilde\eta}=0$ and $t_{\tilde x}+\tilde\eta t_{\tilde u}=0$. 
Then the subsequent split with respect to the new variable $\tilde\eta$ gives the equations $t_{\tilde x}=0$ and $t_{\tilde u}=0$. 
Therefore, $t$ is a function of only $\tilde t$, i.e., $\tilde t$ depends only on $t$, $\tilde t=T(t)$.
Under this condition the function~$F$ expressed in the old variables does not depend on $\eta_t$, i.e., $F_{\eta_t}=0$ and hence 
$(F_{\eta_x},F_{\eta_u})\ne(0,0)$.

Collecting the coefficients of $\eta_{uu}$ and $\eta_{xu}$ in ${\rm DE}_0'$ gives the system 
\[
(\tilde Qu)F_{\eta_u}-\eta^2(\tilde Qx)F_{\eta_x}=0, \quad
(\tilde Qx)F_{\eta_u}+(\tilde Qu)F_{\eta_x}-2\eta(\tilde Qx)F_{\eta_x}=0. 
\]
Since $(F_{\eta_x},F_{\eta_u})\ne(0,0)$, the determinant of the matrix of this system considered as a system 
of linear algebraic equations with respect to $(F_{\eta_x},F_{\eta_u})$ has to vanish, i.e., 
$(\tilde Qu-\eta\tilde Qx)^2=0$ that implies $\tilde Qu=\eta\tilde Qx$.
Assuming $x$ and $u$ to be functions of $(\tilde t,\tilde x,\tilde u,\tilde\eta)$, 
we expand the expression $\tilde Qx$ and $\tilde Qu$ similarly to $\tilde Qt$ and split the equation $\tilde Qu=\eta\tilde Qx$ 
with respect to the new jet variables $\tilde\eta_{\tilde x}$ and $\tilde\eta_{\tilde u}$. 
This results to the equations $u_{\tilde\eta}=\eta x_{\tilde\eta}$ and 
$u_{\tilde x}+\tilde\eta u_{\tilde u}=\eta(x_{\tilde x}+\tilde\eta x_{\tilde u})$. 
Alternating the old and new variables in any derived equation gives a correct equation. 
Therefore, we also have the equations $U_\eta=HX_\eta$, $U_x+\eta U_u=H(X_x+\eta X_u)$. 

The next term for collecting coefficients in ${\rm DE}_0'$ is $\eta_t$. 
The equation obtained by this split is presented as $AG=\tilde A(\tilde Qx)F_{\eta_x}$, 
where $G$ denotes the coefficient of $\eta_t$ in $\tilde\eta_{\tilde t}$. 
Under the above-stated conditions, the expressions appearing in this equation take the form 
\[
F=\frac1{\Delta}\left(\frac{D(H,U)}{D(x,u)}+H\frac{D(X,H)}{D(x,u)}\right),\quad
G=\frac1{T_t\Delta}\frac{\p(H,X,U)}{\p(\eta,x,u)},\quad
\tilde Qx=\frac{U_u-HX_u}{\Delta}. 
\] 
Hereafter $\Delta=D(X,U)/D(x,u)$ $(\ne0)$, and 
\[
\frac{\p(Z^1,\dots,Z^k)}{\p(z_1,\dots,z_k)} \quad\mbox{and}\quad \frac{D(Z^1,\dots,Z^k)}{D(z_1,\dots,z_k)}
\]
denote the usual and total Jacobians of the functions $Z^1$, \dots, $Z^k$ with respect the variables $z_1$, \dots, $z_k$, respectively.  
Note that in the case of a single dependent variable each total Jacobian is, at most, 
a first-order polynomial in derivatives of this dependent variable. 
Removing the denominators from the equations $AG=\tilde A(\tilde Qx)F_{\eta_x}$ results in the equation 
\begin{gather*}
A\Delta^2\dfrac{\p(H,X,U)}{\p(\eta,x,u)}=\tilde A(U_u-HX_u)\left[
\Delta\left(\dfrac{\p(H,U)}{\p(\eta,u)}+H\dfrac{\p(X,H)}{\p(\eta,u)}\right)\right.
\\[1ex]\phantom{A\Delta^2\dfrac{\p(H,X,U)}{\p(\eta,x,u)}=}
\left.{}-\dfrac{\p(X,U)}{\p(\eta,u)}\left(\dfrac{D(H,U)}{D(x,u)}+H\dfrac{D(X,H)}{D(x,u)}\right)\right]
\end{gather*}
the right-hand side of which is at most a first-order polynomial in $\eta_x$ and $\eta_u$.
In view of nonvanishing $A$ and $\p(H,X,U)/\p(\eta,x,u)$, this implies that the coefficients of $\eta_x$ and $\eta_u$ in $\Delta$ equal zero, 
i.e.,  $\p(X,U)/\p(\eta,u)=0$ and $\p(X,U)/\p(x,\eta)=0$. 
Then $\p(X,U)/\p(x,u)\ne0$ since otherwise the transformation $\mathcal T$ is singular.
Hence $X_\eta=U_\eta=0$.

Collecting the coefficients of $\eta_x^2$ in ${\rm DE}_0'$ leads to the equation $H_{\eta\eta}(U_u-HX_x)^2=0$. 
Note that $U_u-HX_x=(\tilde Qx)\Delta\ne0$ since $\Delta\ne0$ and $\tilde Qx\ne0$. 
(Via the split with respect to unconstrained tilde variables, vanishing $\tilde Qx$ implies the condition $x_{\tilde x}=x_{\tilde u}=x_{\tilde\eta}=0$ 
which contradict the nonsingularity of the inverse of $\mathcal T$.)
Therefore, $H_{\eta\eta}=0$, i.e., $H=H^1(t,x,u)\eta+H^0(t,x,u)$, where $H^1=H_\eta\ne0$.
Knowing the explicit dependence of $H$ on $\eta$ allow us to additionally split all equations with respect to $\eta$. 
Thus, splitting the equation $U_x+\eta U_u=H(X_x+\eta X_u)$ gives the condition $X_u=0$ (hence $X_xU_u\ne0$) 
and, then, the conditions $H^1=U_u/X_x$ and $H^0=U_x/X_x$. 
The equation ${\rm DE}_0'$ contains only a single term including $\eta^2\eta_u$. 
Equating the corresponding coefficient to zero, we derive the condition $U_{uu}=0$.  

The whole set of the above found conditions on $T$, $X$, $U$ and $H$ implies that the form of the transformation $\mathcal T$ is described 
by formulas~\eqref{EqGenFormOfTransOfInhomLPEs} and~\eqref{EqTransOfCoeffsOfRedOps0ForLPEs}. 
Then the operator $Q$ is transformed in a simple way: $\tilde Q=X_x^{\,-1}Q$. 
This gives us the idea to rewrite the equations ${\rm DE}_0(\mathcal L)$ and ${\rm DE}_0(\tilde{\mathcal L})$ 
in terms of the operators $Q$ and $\tilde Q$, respectively. 
Thus, the equation ${\rm DE}_0(\mathcal L)$ has the form 
\[
\eta_t+\eta_u(AQ\eta+B\eta+Cu)=AQ^2\eta+(A_x+B)Q\eta+(B_x+C)\eta+C_xu.
\]
All derivatives of $\eta$ containing the differentiation with respect to~$x$ are excluded from ${\rm DE}_0'$ by the substitution $\eta_x=Qu-\eta\eta_u$, 
and hence ${\rm DE}_0'$ can be split with respect to $Q^2\eta$, $\eta_u$, $Q\eta$, $\eta$ and~$u$. 
Collecting the coefficients of the terms $\eta_uQ\eta$, $\eta_u\eta$, $\eta_uu$ and $\eta_u$, 
we obtain the formulas~\eqref{EqTransOfCoeffsOfLPE} for transformations of the arbitrary elements $A$, $B$ and $C$ and the condition $L(U^1/U^0)=0$.
\end{proof}

\begin{note}
We do not split under deriving determining equations in the proof of Lemma~\ref{LemmaOnAdmTransOfDetSystemOfRedOps0ForLPEs} as much as possibly 
since the resulting system would be too cumbersome and, moreover, the proof of the next Theorem~\ref{TheoremOnAdmTransOfDetSystemOfRedOps0ForLPEs}
implies that in fact this complete system is reduced to the set of conditions presented in Lemma~\ref{LemmaOnAdmTransOfDetSystemOfRedOps0ForLPEs}.   
\end{note}

\begin{theorem}\label{TheoremOnAdmTransOfDetSystemOfRedOps0ForLPEs}  
There exists a canonical one-to-one correspondence between the sets of admissible transformations 
of classes~\eqref{EqGenLPE} and~\eqref{EqDet0ForRedOpsOfLPEs}. 
Namely, each point transformation between equations $\mathcal L$ and $\tilde{\mathcal L}$ from class~\eqref{EqGenLPE} 
is extended to a point transformation between the associated equations ${\rm DE}_0(\mathcal L)$ and ${\rm DE}_0(\tilde{\mathcal L})$ 
according to formula~\eqref{EqTransOfCoeffsOfRedOps0ForLPEs}. 
In both the transformations the variables $(t,x,u)$ and the arbitrary elements are transformed in the same way. 
The extended transformations exhaust the sets of admissible transformation in class~\eqref{EqDet0ForRedOpsOfLPEs}.
\end{theorem} 

\begin{proof}
The extension of each admissible transformation in class~\eqref{EqGenLPE} by formula~\eqref{EqTransOfCoeffsOfRedOps0ForLPEs} 
gives an admissible transformation in class~\eqref{EqDet0ForRedOpsOfLPEs}. 
Therefore, it is enough to check that every admissible transformation in class~\eqref{EqDet0ForRedOpsOfLPEs} 
coincides with the extension of an admissible transformation in class~\eqref{EqGenLPE}. 
We take two equations from class~\eqref{EqDet1ForRedOpsOfLPEs} which are connected via a point transformation. 
They necessarily are determining equations for reduction operators with the zero coefficients of $\p_t$ and the unit coefficients of $\p_x$
for some equations $\mathcal L$ and $\tilde{\mathcal L}$ from class~\eqref{EqGenLPE}. 
Therefore, these equations can be denoted by ${\rm DE}_0(\mathcal L)$ and ${\rm DE}_0(\tilde{\mathcal L})$, respectively.
Consider a point transformation~$\breve{\mathcal T}$ mapping ${\rm DE}_0(\mathcal L)$ to ${\rm DE}_0(\tilde{\mathcal L})$. 
In view of Lemma~\ref{LemmaOnAdmTransOfDetSystemOfRedOps0ForLPEs}, the transformation~$\breve{\mathcal T}$ 
has the form given by formulas~\eqref{EqGenFormOfTransOfInhomLPEs} and~\eqref{EqTransOfCoeffsOfRedOps0ForLPEs} 
and, therefore, is projectable on the space of the variables $(t,x,u)$. 
Denote its projection by~$\mathcal T$.
The pairs $({\rm DE}_0(\mathcal L),{\rm DE}_0(\tilde{\mathcal L}))$ and $(\mathcal L,\tilde{\mathcal L})$ have the same tuples of 
arbitrary elements transformed by the same formulas~\eqref{EqTransOfCoeffsOfLPE}. 
Hence Lemmas~\ref{LemmaOnAdmTransOfInhomLPEs} and~\ref{LemmaOnAdmTransOfDetSystemOfRedOps0ForLPEs} imply 
that $\mathcal T$ is a point transformation from $\mathcal L$ to $\tilde{\mathcal L}$. 
It is clear that the transformation~$\breve{\mathcal T}$ is the extension of~$\mathcal T$ by formula~\eqref{EqTransOfCoeffsOfRedOps0ForLPEs}, 
i.e., $\breve{\mathcal T}=\mathcal T_0$.  
\end{proof}

\begin{corollary}\label{CorollaryOnEquivGroupOfDetSystemOfRedOps0ForLPEs}
The equivalence group $G^\sim_0$ of class~\eqref{EqDet0ForRedOpsOfLPEs} is isomorphic to the equivalence group $G^\sim$ of class~\eqref{EqGenLPE}. 
The canonical isomorphism is established by the extension of transformations from $G^\sim_0$ to the variable~$\eta$
via formula~\eqref{EqTransOfCoeffsOfRedOps0ForLPEs}, where $U^0=0$. 
\end{corollary}

\begin{corollary}\label{CorollaryOnLieSymAlgOfDetSystemOfRedOps0ForLPEs}
For any equation~$\mathcal L$ from class~\eqref{EqGenLPE}, the maximal point symmetry groups 
(resp. the maximal Lie invariance algebras) of the equations~$\mathcal L$ and ${\rm DE}_0(\mathcal L)$ are isomorphic. 
The canonical isomorphism between the algebras is realized via the extension of 
each Lie symmetry operator $Q=\tau\p_t+\xi\p_x+(\zeta^1u+\zeta^0)\p_u$ of~$\mathcal L$ 
to the Lie symmetry operator $Q_1=Q+\bigl((\zeta^1-\xi_x)\eta+\zeta^1_xu+\zeta^0_x\bigr)\p_\eta$ of ${\rm DE}_0(\mathcal L)$. 
\end{corollary}

In view of Corollaries~\ref{CorollaryOnEquivGroupOfDetSystemOfRedOps0ForLPEs} and \ref{CorollaryOnLieSymAlgOfDetSystemOfRedOps0ForLPEs}, 
the results on the group classification of class~\eqref{EqDet0ForRedOpsOfLPEs} follow from Theorem~\ref{TheoremOnGroupClassificationOfLPEs}.
  
\begin{corollary}\label{CorollaryOnGroupClassificationOfDetSystemOfRedOps0ForLPEs}
The kernel Lie algebra of class~\eqref{EqDet0ForRedOpsOfLPEs} is $\langle I_0\rangle$, where $I_0=u\p_u+\eta\p_\eta$.
Any equation from class~\eqref{EqDet0ForRedOpsOfLPEs} is invariant with respect to the operators of the form $Z_0(f)=f\p_u+f_x\p_\eta$,
where the parameter-function $f=f(t,x)$ runs through the solution set of the associated equation 
$f_t=Af_{xx}+Bf_x+Cf$.
All possible $G^\sim_0$-inequivalent cases of extension of the maximal 
Lie invariance algebra are exhausted by the following equations of the reduced form~\eqref{EqReducedDet0ForRedOpsOfLPEs} 
(the values of~$V$ are given together with the corresponding maximal Lie invariance algebras):

\vspace{1.5ex}

$\makebox[6mm][l]{\rm 1.}
V=V(x) \colon\quad \langle\p_t,\,I_0,\,Z_0(f)\rangle;$

\vspace{1.5ex}

$\makebox[6mm][l]{\rm 2.}
V=\mu x^{-2},\ \mu\ne0 \colon\quad \langle\p_t,\, D_0,\, \Pi_0,\, I_0,\,Z_0(f)\rangle;$

\vspace{1.5ex}

$\makebox[6mm][l]{\rm 3.}
V=0 \colon\quad \langle\p_t,\, \p_x,\, G_0,\, D_0,\, \Pi_0,\, I_0,\,Z_0(f)\rangle$.

\vspace{1.5ex}

\noindent
Here 
\begin{gather*}
D_0=2t\p_t+x\p_x-\eta\p_\eta, 
\\
\Pi_0=4t^2\p_t+4tx\p_x-(x^2+2t)u\p_u-(x\eta+6t\eta+2xu)\p_\eta, 
\\
G_0=2t\p_x-xu\p_u-(x\eta+u)\p_\eta.
\end{gather*}
\end{corollary}

\section{Lie reductions of determining equations}\label{SectionOnLieReductionsOfDetEqsOfReditionOpsOfLPEs}

Suppose that an equation~$\mathcal L$ from class~\eqref{EqGenLPE} admits a Lie symmetry operator
$Q=\tau\p_t+\xi\p_x+\zeta\p_u$. The coefficients of~$Q$ necessarily satisfy the conditions
$\tau_x=\tau_u=0$, $\xi_u=0$ and $\zeta_{uu}=0$, i.e., 
$\tau=\tau(t)$, $\xi=\xi(t,x)$ and $\zeta=\zeta^1(t,x)u+\zeta^0(t,x)$, 
and $\zeta^0$ is a solution of~$\mathcal L$.

In view of Corollaries~\ref{CorollaryOnLieSymAlgOfDetSystemOfRedOps1ForLPEs} and \ref{CorollaryOnLieSymAlgOfDetSystemOfRedOps0ForLPEs}, 
the determining equations ${\rm DE}_1(\mathcal L)$ and ${\rm DE}_0(\mathcal L)$ respectively possess 
the Lie symmetry operators $Q_1$ and $Q_0$ associated with~$Q$, 
which can be applied to reduce the determining equations and construct their exact solutions. 
The found solutions of the determining equations give the reduction operators of a special kind for the initial equation~$\mathcal L$, 
implicitly connected with Lie invariance properties of~$\mathcal L$. 
The question is what properties the solutions of~$\mathcal L$, invariant with respect to such reduction operators, possess, e.g., 
whether these solutions necessarily are Lie invariant or they are not. 

An admissible transformation~$\mathcal T$ of the equation~$\mathcal L$ in class~\eqref{EqGenLPE} 
has form~\eqref{EqGenFormOfTransOfInhomLPEs} and maps the pair $(\mathcal L,Q)$ to a pair $(\mathcal L',Q')$, 
where the equation~$\mathcal L'$ also belongs to class~\eqref{EqGenLPE} 
and $Q'$ is a nontrivial (resp. trivial) Lie symmetry operator of~$\mathcal L'$ if  
$Q$ is a nontrivial (resp. trivial) Lie symmetry operator of~$\mathcal L$. 
Up to the equivalence generated by the set of all admissible transformations of class~\eqref{EqGenLPE} 
(see Lemma~\ref{LemmaOnAdmTransOfInhomLPEs}) in the set of pairs `(equation of form~\eqref{EqGenLPE}, its Lie symmetry operator)', 
we can assume that $Q\in\{\p_t,\p_x\}$ or $Q\in\{u\p_u,\p_u\}$ if $Q$ is a nontrivial or trivial Lie symmetry operator of~$\mathcal L$, respectively. 
$Q\sim\p_t$ if $\tau\ne0$ and $Q\sim\p_x$ if $\tau=0$ and $\xi\ne0$.

If $Q\in\{\p_t,\p_x\}$, 
the Lie symmetry operator $Q_1$ of the system ${\rm DE}_1(\mathcal L)$ and  
the Lie symmetry operator $Q_0$ of the equation ${\rm DE}_0(\mathcal L)$, which are associated with the operator~$Q$, 
formally have the same form as the operator~$Q$ but are defined in different spaces of variables.

\begin{proposition}
Suppose that an equation~$\mathcal L$ from class~\eqref{EqGenLPE} possesses a Lie symmetry operator $Q=\tau\p_t+\xi\p_x+\zeta\p_u$, 
where necessarily $\tau=\tau(t)$, $\xi=\xi(t,x)$ and $\zeta=\zeta^1(t,x)u+\zeta^0(t,x)$ and additionally $\tau\ne0$. 
Let~$Q_1$ be the associated Lie symmetry operator of the system ${\rm DE}_1(\mathcal L)$, 
a solution $(g^1,g^2,g^3)$ of ${\rm DE}_1(\mathcal L)$ be $Q_1$-invariant and 
$R=\p_t+g^1\p_x+(g^2u+g^3)\p_u\in\mathcal Q_1(\mathcal L)$ be the corresponding reduction operator. 
Then the functions~$g^1$, $g^2$ and~$g^3$ are expressed, according to formulas~\eqref{EqTrans1BetweenInitialAndDetEqsForLPEs}, 
via a solution $(v^1,v^2,v^3)$ of the uncoupled system $3\mathcal L$, 
which is invariant with respect to the Lie symmetry operator 
\[
\tau\p_t+\xi\p_x+\zeta^1v^1\p_{v^1}+\zeta^1v^2\p_{v^2}+(\zeta^1v^3+\zeta^0)\p_{v^3}+\lambda_{iq}v^q\p_{v^i}
\]
of this system for some constants $\lambda_{iq}$, $i=1,2,3$, $q=1,2$. 
Here the functions~$v^1$ and $v^2$ have to be linearly independent.
Each $R$-invariant solution of~$\mathcal L$ is a linear combination, with the unit coefficient of~$v^3$, 
of the components of the Lie invariant solution~$(v^1,v^2,v^3)$ of the system $3\mathcal L$.
\end{proposition}

\begin{proof}
It is sufficient to consider only the reduced form of Lie symmetry operators, 
which is $Q=\p_t$ in the case $\tau\ne0$.
Then $Q_1=\p_t$.
The equation~$\mathcal L$ is Lie invariant with respect to the operator~$\p_t$ if and only if $A_t=B_t=C_t=0$. 
Consider an operator $R=\p_t+g^1\p_x+(g^2u+g^3)\p_u\in\mathcal Q_1(\mathcal L)$, where 
the coefficient tuple $(g^1,g^2,g^3)$ is a $Q_1$-invariant solution of ${\rm DE}_1(\mathcal L)$, i.e., 
it additionally satisfies the condition $g^1_t=g^2_t=g^3_t=0$. 
An ansatz constructed with the operator~$R$ has the form $u=f^1(x)\varphi(\omega)+f^0(x)$, 
where $f^1=f^1(x)\ne0$ and $f^0=f^0(x)$ are given coefficients, $\varphi=\varphi(\omega)$ is the new unknown function, 
$\omega=t+\varrho(x)$ is the invariant independent variable and $\varrho_x\ne0$.
This ansatz reduces $\mathcal L$ to an (in general, inhomogeneous) linear second-order constant-coefficient 
ordinary differential equation in~$\varphi$, which we denote by $\mathcal L'$. 
The general solution of~$\mathcal L'$ is represented in the form 
$\varphi=c_p\varphi^p(\omega)+\varphi^3(\omega)$, 
where $\varphi^3$ is a particular solution of $\mathcal L'$,  
$\varphi^1$ and $\varphi^2$ are linearly independent solutions of the corresponding homogeneous equation 
and $c_1$ and $c_2$ are arbitrary constants. 
Let us recall that $p,q=1,2$.
Substituting the general solution of $\mathcal L'$ into the ansatz, 
we obtain the two-parametric family of solutions of~$\mathcal L$, having form~\eqref{EqLinear2ParamSetOfSolutionsOfLPEs} 
with $v^p=f\varphi^p$ and $v^3=f\varphi^3+g$. 
Due to $\mathcal L'$ is a constant-coefficient equation, the functions $v^i$ admit 
the representation $v^p=\psi^{pq}(t)\theta^q(x)$ and $v^3=\psi^{3q}(t)\theta^q(x)+\theta^3(x)$, where 
$\psi^{iq}_t=\lambda_{ip}\psi^{pq}$ for some constants $\lambda_{ip}$ depending on the coefficients of~$\mathcal L'$. 
Therefore $(v^1,v^2,v^3)$ is a solution of the system $3\mathcal L$, 
which is invariant with respect to the Lie symmetry operator $\p_t+\lambda_{iq}v^q\p_{v^i}$ of this system.
\end{proof}

\begin{proposition}
Suppose that the system ${\rm DE}_1(\mathcal L)$ associated with an equation~$\mathcal L$ from class~\eqref{EqGenLPE}
possesses a Lie invariance operator~$Q_1$ with the vanishing coefficient of~$\p_t$ and a nonvanishing coefficient of~$\p_x$. 
Let a solution $(g^1,g^2,g^3)$ of ${\rm DE}_1(\mathcal L)$ be invariant with respect to~$Q_1$. 
Then the associated reduction operator $\p_t+g^1\p_x+(g^2u+g^3)\p_u$ of the equation~$\mathcal L$ 
is necessarily equivalent to a Lie invariance operator of~$\mathcal L$.
\end{proposition}

\begin{proof}
Consider the case $Q=\p_x$.
The equation~$\mathcal L$ possesses the Lie symmetry operator~$\p_x$ if and only if $A_x=B_x=C_x=0$. 
Then the equivalence transformation $\tilde t=T(t)$, $\tilde x=x+\varphi(t)$ and $\tilde u=\psi(t)u$, 
where $T_t=A$, $\varphi_t=B$ $\psi_t=C\psi$ and $\psi\ne0$, maps $Q$ to $\p_{\tilde x}$ and 
reduces $\mathcal L$ to the linear heat equation 
$\tilde u_{\tilde t}=\tilde u_{\tilde x\tilde x}$ associated with the values $\tilde A=1$ and $\tilde B=\tilde C=0$. 
That is why without loss of generality we can assume that $A=1$ and $B=C=0$. 
An ansatz constructed for the system ${\rm DE}_1(\mathcal L)$ by the operator~$Q_1=\p_x$ is 
$g^i=g^i(t)$ and the corresponding reduced system has the form $g^i_t=0$, i.e., $g^i=\const$.
The operator $\p_t+g^1\p_x+(g^2u+g^3)\p_u$ with constant coefficients 
belongs to the maximal Lie invariance algebra of the equation~$\mathcal L$ 
which coincides under our suppositions with the linear heat equation.  
The obtained statement is reformulated for the general form of~$Q$ with the vanishing coefficient of~$\p_t$. 
\end{proof}

Results on Lie solutions of the determining equation ${\rm DE}_0(\mathcal L)$
can be presented as a single statement 
without split into different cases depending on the structure of the corresponding Lie symmetry operators. 
To formulate them in a compact form, we need to introduce at first the auxiliary notion of 
\emph{one-parametric solution families} of the equation~$\mathcal L$, 
\emph{associated with the Lie symmetry operator}~$Q$ of~$\mathcal L$. 
The set of such families is partitioned into two subsets which are respectively formed by 
the \emph{singular associated families} consisting of $Q$-invariant solutions of~$\mathcal L$ and 
the \emph{regular associated families} obtained via acting on fixed non-$Q$-invariant solutions of~$\mathcal L$ 
by the one-parametric transformation group generated by~$Q$.

Let us recall that $Q_0$ denotes the Lie symmetry operator of ${\rm DE}_0(\mathcal L)$, associated with~$Q$. 
Equivalent families of solutions, which differ only by parametrization, are identified. 
In particular, regular one-parametric families associated with the same operator are equivalent if and only if 
they differ only by parameter shifts. 
Such families are obtained by the action of the same one-parametric transformation group on fixed solutions which are similar 
with respect to this group. 
A neighborhood of a nonsingular point of~$Q$ is considered. 
(Otherwise, the one-to-one correspondence in the next theorem may be broken. 
In some cases it can saved by taking into account discrete symmetry transformations, 
see Note~14 of~\cite{Popovych2007a}.)

Formulas~\eqref{EqTransOfCoeffsOfRedOps0ForLPEs} and~\eqref{EqCoeffsOfSymOpersOfDetSystemOfRedOps0ForLPEs} imply 
the following statement which will be used below.

\begin{proposition}\label{PropositionOnInverianceOfAdditionalConstraint0AndEtaForLPEs}
Let an equation~$\mathcal L$ from class~\eqref{EqGenLPE} be invariant with respect to a point transformation~$\mathcal T$ 
(resp. an operator~$Q$) and  
the function $\eta=\eta(t,x,u)$ be a solution of the associated determining equation ${\rm DE}_0(\mathcal L)$. 
Then the equations $u_x=\eta(t,x,u)$ admits the transformation~$\mathcal T$ (resp. the operator~$Q$) 
as a point symmetry transformation (resp. a Lie symmetry operator) if and only if 
the function $\eta$ is an invariant of the associated point symmetry transformation~$\mathcal T_0$ 
(resp. the associated Lie symmetry operator~$Q_0$) of the equation ${\rm DE}_0(\mathcal L)$.
\end{proposition}

\begin{theorem}\label{TheoremOnRedOps0OfLPEsInvWrtLieOps}
For each equation~$\mathcal L$ from class~\eqref{EqGenLPE} and each Lie symmetry operator~$Q$ of~$\mathcal L$, 
there exists a one-to-one correspondence between $Q_0$-invariant solutions of the determining equation ${\rm DE}_0(\mathcal L)$ 
and one-parametric families of solutions of~$\mathcal L$, associated with~$Q$.
Namely, the reduction of the equation~$\mathcal L$ by an operator $\p_x+\eta\p_u$, 
where the coefficient~$\eta$ is a $Q_0$-invariant solution of ${\rm DE}_0(\mathcal L)$, 
gives a one-parametric solution family of~$\mathcal L$, associated with~$Q$. 
And vice versa, each family of the above kind consists of solutions invariant with respect to an operator $\p_x+\eta\p_u$, 
where the coefficient~$\eta$ is a $Q_0$-invariant solution of ${\rm DE}_0(\mathcal L)$.
\end{theorem}

\begin{proof}
Suppose that an equation~$\mathcal L$ from class~\eqref{EqGenLPE} admits a Lie symmetry operator~$Q$. 
We denote the one-parametric transformation group with~the infinitesimal operator~$Q$ by $G$.
Let a solution~$\eta$ of the equation ${\rm DE}_0(\mathcal L)$ be invariant with respect to the 
associated operator $Q_0$. 
Then the system~$\mathcal L_\eta$ of the equation~$\mathcal L$ with the additional constraint $u_x=\eta$ 
possesses $Q$ as a Lie symmetry operator. 
The general solution~$\mathcal F$ of~$\mathcal L_\eta$ is a one-parametric solution family of~$\mathcal L$.   
There are two different cases of the structure of~$\mathcal F$. 
In the first case the family~$\mathcal F$ consists of $Q$-invariant solutions of~$\mathcal L$ and, therefore,  
is a singular one-parametric solution family associated with the operator~$Q$. 
In the second case the family~$\mathcal F$ contains a solution $u=u^0(t,x)$ of~$\mathcal L$, which is not $Q$-invariant. 
A one-parametric family of solutions of~$\mathcal L_\eta$ obtained via acting on the solution $u^0$ by transformations from~$G$ 
is equivalent to~$\mathcal F$. 
Therefore, $\mathcal F$ is a regular one-parametric solution family associated with the operator~$Q$.

Vice versa, if a one-parametric solution family of the equation~$\mathcal L$ is associated with the operator~$Q$ then 
the corresponding additional constraint $u_x=\eta$ with a solution~$\eta$ of ${\rm DE}_0(\mathcal L)$ 
admits $Q$ as a Lie symmetry operator. 
In view of Proposition~\ref{PropositionOnInverianceOfAdditionalConstraint0AndEtaForLPEs},
this implies that the function $\eta$ is $Q$-invariant.
\end{proof}

Since the determining equation ${\rm DE}_0(\mathcal L)$ has three independent variable, 
it also admits Lie reductions with respect to two-dimensional subalgebras of its maximal Lie invariance algebras 
to ordinary differential equations and, therefore, possesses the corresponding invariant solutions. 
To formulate the statement on such solutions analogously to Theorem~\ref{TheoremOnRedOps0OfLPEsInvWrtLieOps},
we need to define \emph{one-parametric families of solutions} of the equation~$\mathcal L$, 
\emph{associated with the two-dimensional Lie invariance algebra}~$\mathfrak g$ of~$\mathcal L$. 
The whole set of associated families is also partitioned into the subsets of the \emph{singular} and \emph{regular} families. 
Each singular associated family consists of $\mathfrak g$-invariant solutions of~$\mathcal L$. 
Each regular associated family is obtained via acting on fixed $Q^1$-invariant and non-$Q^2$-invariant solution of~$\mathcal L$ 
by the one-parametric transformation group generated by~$Q^2$. 
Here $Q^1$ and $Q^2$ are arbitrary linearly independent elements of~$\mathfrak g$.

\begin{theorem}\label{TheoremOnRedOps0OfLPEsInvWrt2DLieInvAlgs}
Suppose that a two-dimensional Lie invariance algebra $\mathfrak g$ of an equation~$\mathcal L$ from class~\eqref{EqGenLPE} 
induces the Lie invariance algebra $\mathfrak g_0$ of the corresponding determining equation ${\rm DE}_0(\mathcal L)$, 
which is appropriate for Lie reduction of ${\rm DE}_0(\mathcal L)$. 
Then there exists a one-to-one correspondence between $\mathfrak g_0$-invariant solutions of ${\rm DE}_0(\mathcal L)$ 
and one-parametric families of solutions of~$\mathcal L$, associated with~$\mathfrak g$.
Namely, the reduction of~$\mathcal L$ by an operator $\p_x+\eta\p_u$, 
where the coefficient~$\eta$ is a $\mathfrak g_0$-invariant solution of ${\rm DE}_0(\mathcal L)$, 
gives a one-parametric family of solutions of~$\mathcal L$, associated with~$\mathfrak g$. 
And vice versa, each family of this kind consists of solutions invariant with respect to an operator $\p_x+\eta\p_u$, 
where the coefficient~$\eta$ is a $\mathfrak g_0$-invariant solution of ${\rm DE}_0(\mathcal L)$.
\end{theorem}

\begin{proof}
We denote by $G$ the two-parametric transformation group with the Lie algebra~$\mathfrak g$ 
and locally parameterize elements of $G$ in a neighborhood of the identical transformation by the pair $(\varepsilon_1,\varepsilon_2)$: 
$g(\varepsilon_1,\varepsilon_2)\in G$. In particular, $g(0,0)$ is the identical transformation 
and the infinitesimal operators $Q^i=g_{\varepsilon_i}(0,0)$, $i=1,2$, form a basis of the algebra~$\mathfrak g$.
Let a solution~$\eta$ of the equation ${\rm DE}_0(\mathcal L)$ be invariant with respect to the 
associated algebra $\mathfrak g_0$. 
Then $\mathfrak g$ is a Lie invariance algebra of 
the system~$\mathcal L_\eta$ formed by the equation~$\mathcal L$ and the additional constraint $u_x=\eta$. 
The general solution~$\mathcal F$ of~$\mathcal L_\eta$ is a one-parametric solution family of~$\mathcal L$. 
We explicitly represent this family by the formula $u=f(t,x,\varkappa)$. 
There are two different cases of its possible structure. 
The family~$\mathcal F$ can consist of $\mathfrak g$-invariant solutions of~$\mathcal L$ and, therefore,  
be a singular one-parametric solution family associated with the algebra~$\mathfrak g$. 
The other possibility is that the family~$\mathcal F$ contains a solution $u=f(t,x,\varkappa_0)$ of~$\mathcal L$, which is not $\mathfrak g$-invariant. 
Then the solution $u=f(t,x,\varkappa_0)$ is invariant with respect to the operator $\varkappa_{0,1}Q^2-\varkappa_{0,2}Q^1\in\mathfrak g$, 
where $\varkappa_{0,i}=(g(\varepsilon_1,\varepsilon_2)\varkappa_0)_{\varepsilon_i}|_{(\varepsilon_1,\varepsilon_2)=(0,0)}$.
The action of the one-parametric subgroup~$G'$ of~$G$ with the infinitesimal operator~$\varkappa_{0,1}Q^1+\varkappa_{0,2}Q^2\in\mathfrak g$ 
is (locally) transitive on~$\mathcal F$. 
It means that $\mathcal F$ is a regular one-parametric solution family associated with the algebra~$\mathfrak g$, 
which is obtained via acting by~$G'$ on the fixed solution $u=f(t,x,\varkappa_0)$.

Conversely, if a one-parametric solution family of the equation~$\mathcal L$ is associated with the algebra~$\mathfrak g$ then 
the corresponding additional constraint $u_x=\eta$, where $\eta=\eta(t,x,u)$ is a solution of ${\rm DE}_0(\mathcal L)$, 
admits $\mathfrak g$ as a Lie symmetry algebra. 
In view of Proposition~\ref{PropositionOnInverianceOfAdditionalConstraint0AndEtaForLPEs},
this implies that the function $\eta$ is $\mathfrak g$-invariant.
\end{proof}

\section{Particular cases of reductions and linearization}\label{SectionOnParticularCasesOfReductionsAndLinearization}

All possible Lie reductions of the determining equations~\eqref{EqDet1ForRedOpsOfLPEs} and~\eqref{EqDet0ForRedOpsOfLPEs} 
for reduction operators and the corresponding invariant solutions of equations from class~\eqref{EqGenLPE} 
are studied in the previous section. 
Now we consider a few examples of typical additional conditions to the determining equations, which are 
different from Lie ones. 
A special attention is paid to an interpretation of the confinement of the linearizing transformations given in 
Corollaries~\ref{CorollaryOnLinearizationOfDetEqs1ForRedOpsOfLPEs} and~\ref{CorollaryOnLinearizationOfDetEqs0ForRedOpsOfLPEs} 
to the particular cases under consideration. 
Presented examples also show that nontrivial reduction operators associated with nontrivial additional conditions 
to the determining equations can finally leads to trivial solutions of equations from class~\eqref{EqGenLPE}.

We fix an equation~$\mathcal L$ from class~\eqref{EqGenLPE}.
The extension of possibilities for constraints of the determining equations in comparison with 
the initial equation~$\mathcal L$ is connected with a greater number of unknown functions in ${\rm DE}_1(\mathcal L)$ 
and the additional independent variable~$u$ in ${\rm DE}_0(\mathcal L)$.

Consider at first reduction operators of~$\mathcal L$ with the vanishing coefficients of~$\p_t$.

\begin{example}
Suppose that $Q_0=\p_x$ is a reduction operator of ${\rm DE}_0(\mathcal L)$. It means that 
the arbitrary elements satisfy the condition $A_x=B_{xx}=C_x=0$. 
The problem is to investigate solutions of ${\rm DE}_0(\mathcal L)$, which are invariant with respect to~$Q_0$. 
We do an equivalence transformation of the form 
\[\tilde t=T(t),\quad \tilde x=X^1(t)x+X^0(t),\quad \tilde u=U^1(t)u,\] 
where the arbitrary elements~$A$, $B$ and~$C$ and the function~$\eta$ are transformed 
according to formulas~\eqref{EqTransOfCoeffsOfLPE} and~\eqref{EqTransOfCoeffsOfRedOps0ForLPEs}. 
The parameter-functions $T$, $X^1$, $X^0$ and $U^1$ can be chosen in such a way that 
$\tilde A=1$, $\tilde B=0$ and  $\tilde C=0$. 
In the new variables the operator~$Q_0$ equals $X^1\p_{\tilde x}$ and hence is equivalent to~$\p_{\tilde x}$. 
This is why we can assume without loss of generality that $A=1$, $B=0$ and  $C=0$, 
i.e., $\mathcal L$ coincides with the linear heat equation.
Then $Q_0=\p_x$ is a Lie symmetry operator of ${\rm DE}_0(\mathcal L)$. 
The corresponding reduced equation $\eta_t=\eta\eta_{uu}$ for the function $\eta=\eta(t,u)$ is equivalent, 
on the subset of nonvanishing solutions, to the remarkable nonlinear diffusion equation $\zeta_t=(\zeta^{-2}\zeta_u)_u$, 
where $\zeta=1/\eta$. 
It is well known that this diffusion equation is linearized 
to the linear heat equation~\cite{Bluman&Kumei1980,Storm1951}. 
We derive this transformation via confining the transformation of ${\rm DE}_0(\mathcal L)$ to, formally, $\mathcal L$, 
presented in Corollary~\ref{CorollaryOnLinearizationOfDetEqs0ForRedOpsOfLPEs}. 
We put $\Phi=\Psi(t,u)-x$, where $\Psi_u\ne0$. 
Then 
\[\eta=-\frac{\Phi_x}{\Phi_u}=\frac1{\Psi_u},\]
i.e., $\zeta=\Psi_u$. 
After integrating, we obtain the equation $\Psi_t=\Psi_{uu}/\Psi_u{}^2+\beta(t)$ in the function $\Psi=\Psi(t,u)$. 
The ``integration constant'' $\beta=\beta(t)$ can be assumed to vanish due to the ambiguity in the connection between 
$\zeta$ and $\Psi$. 
The confinement of transformation~\eqref{EqTrans0BetweenInitialAndDetEqs} is the hodograph transformation 
\begin{gather*}
\mbox{the new independent variables:}\qquad\tilde t=t, \quad \tilde x=\Psi, 
\\ 
\lefteqn{\mbox{the new dependent variable:}}\phantom{\mbox{the new independent variables:}\qquad }\tilde u=u 
\end{gather*}
since here the variable~$x$ has to be replaced by $\Psi=x+\Phi$. 
The application of this transformation results in the linear heat equation 
$\tilde u_{\tilde t}=\tilde u_{\tilde x\tilde x}$.
Note that the above interpretation of the confinement of transformation~\eqref{EqTrans0BetweenInitialAndDetEqs}
differs from the interpretation in the proof of Theorem~9 of~\cite{Popovych2007a}.
\end{example}

\begin{example}\label{ExampleOnRedOps0EtaUU0OfLPEs}
Let the function~$\eta$ satisfy the additional condition $\eta_{uu}=0$, i.e., $\eta=\eta^1(t,x)u+\eta^0(t,x)$.
Then the equation ${\rm DE}_0(\mathcal L)$ is reduced to the system 
\begin{equation}\arraycolsep=0ex\label{EqSystemBeingEquivDE0OfLPEsWithEta_uu0}
\begin{array}{l}
\eta^1_t=\bigl(A\eta^1_x+A(\eta^1)^2+B\eta^1+C\bigr)_x,\\[1ex]
\eta^0_t=A(\eta^0_{xx}+2\eta^0\eta^1_x)+A_x(\eta^0_x+\eta^0\eta^1)+(B\eta^0)_x+C\eta^0.
\end{array}
\end{equation}
Putting $\Phi=\Phi^1(t,x)u+\Phi^0(t,x)$, we rewrite the transformation described 
in Corollary~\ref{CorollaryOnLinearizationOfDetEqs0ForRedOpsOfLPEs} in terms of~$\eta^1$ and~$\eta^0$. 
The condition $\eta=-\Phi_x/\Phi_u$ implies that $\eta^1=-\Phi^1_x/\Phi^1$ and $\eta^0=-\Phi^0_x/\Phi^1$. 
The hodograph transformation~\eqref{EqTrans0BetweenInitialAndDetEqs} is equivalent 
to expressing $u$ from the formula for~$\Phi$: 
\[
u=\frac{\Phi-\Phi^0}{\Phi^1}=\Psi^1(t,x)\varkappa+\Psi^0(t,x),
\]
where $\Psi^1=1/\Phi^1$ and $\Psi^0=\Phi^0/\Phi^1$. 
Since the expression for $u$ has to be the solution family of~$\mathcal L$ with the parameter $\varkappa=\Phi$, 
$\Psi^1$ and $\Psi^0$ are solutions of~$\mathcal L$, $\Psi^1\ne0$.  
Finally we derive the representation 
\begin{equation}\label{EqLinearizationTransForDE0OfLPEsWithEta_uu0}
\eta^1=\frac{\Psi^1_x}{\Psi^1},\quad \eta^0=\Psi^0_x-\frac{\Psi^1_x}{\Psi^1}\Psi^0,
\end{equation}
where $\Psi^1$ and $\Psi^0$ are solutions of the initial equation~$\mathcal L$. 
In other words, transformation~\eqref{EqLinearizationTransForDE0OfLPEsWithEta_uu0} reduces 
the nonlinear system~\eqref{EqSystemBeingEquivDE0OfLPEsWithEta_uu0} 
in~$\eta^1$ and~$\eta^0$ to the system of two uncoupled copies of~$\mathcal L$. 
The expression for~$\eta^1$ in~\eqref{EqLinearizationTransForDE0OfLPEsWithEta_uu0} coincides, 
up to sign, with the well-known Cole--Hopf substitution linearizing the Burgers equation.  
(If $A=1$ and $B=C=0$, the first equation of~\eqref{EqSystemBeingEquivDE0OfLPEsWithEta_uu0} 
coincides, up to signs, with the Burgers equation.)
The expression for~$\eta^1$ in~\eqref{EqLinearizationTransForDE0OfLPEsWithEta_uu0} is obtained 
as the Darboux transformation of the solution~$\Psi^0$, associated with the solution~$\Psi^1$. 
It follows from~\eqref{EqLinearizationTransForDE0OfLPEsWithEta_uu0} that 
the reduction operator $R=\p_x+(\eta^1u+\eta^0)\p_u$ is $G^\infty(\mathcal L)$-equivalent 
to the operator $\p_x+\eta^1u\p_u$. 
Indeed, the transformation $\tilde t=t$, $\tilde x=x$, $\tilde u=u-\Psi^0$ belongs to
$G^\infty(\mathcal L)$ and maps the operator~$R$ to 
the operator $\tilde R=\p_{\tilde x}+\eta^1\tilde u\p_{\tilde u}$.

An ansatz constructed with $R$ has the form $u=\Psi^1(t,x)\varphi(\omega)+\Psi^0(t,x)$, 
where $\varphi=\varphi(\omega)$ is an invariant unknown function of the invariant independent variable $\omega=t$. 
The associated reduced equation is $\varphi_\omega=0$, i.e., $\varphi=\const$. 
Therefore, $u=\Psi^1\varkappa+\Psi^0$ is the family of $R$-invariant solutions of~$\mathcal L$. 

Vice versa, the solution family $u=\Psi^1(t,x)\varkappa+\Psi^0(t,x)$ of the equation~$\mathcal L$ is 
necessarily invariant with respect to the reduction operator $\p_x+(\eta^1(t,x)u+\eta^0(t,x))\p_u$, where 
the coefficients~$\eta^1$ and~$\eta^0$ are determined by 
the formulas~\eqref{EqLinearizationTransForDE0OfLPEsWithEta_uu0}. 

As a result, we obtain the following statement.  

\begin{proposition}
For any equation of form~\eqref{EqGenLPE}, 
there exists a one-to-one correspondence between one-parametric families of its solutions, 
linearly depending on parameters, and reduction operators of the form $\p_x+\eta(t,x,u)\p_u$, 
where $\eta_{uu}=0$.
Namely, each operator of such kind corresponds to
the family of solutions which are invariant with respect to this operator. 
\end{proposition}

\end{example}

\begin{example}
At first sight, the additional condition $\eta_x+\eta\eta_u=0$ seems much more complicated than the conditions 
studied in the previous examples. 
In fact, it leads only to solutions of the initial equation~$\mathcal L$, which are first-order polynomials 
with respect to~$x$. 
To see this, we carry out the transformation described 
in Corollary~\ref{CorollaryOnLinearizationOfDetEqs0ForRedOpsOfLPEs} 
and, as a result, obtain the condition $\tilde u_{\tilde x\tilde x}=0$. 
In contrast to the solutions of~$\mathcal L$, the associated solutions of ${\rm DE}_0(\mathcal L)$ 
have a complex structure and are difficult to construct.

The system $S$ consisting of the equations ${\rm DE}_0(\mathcal L)$ and $\eta_x+\eta\eta_u=0$ 
has the compatibility condition $(B_{xx}+2C_x)\eta+C_{xx}u=0$.
Before considering the possible cases, we note that the equation $\eta_x+\eta\eta_u=0$ is invariant with respect to the transformations 
from the equivalence group $G^\sim_0$ of class~\eqref{EqDet0ForRedOpsOfLPEs}, which additionally satisfy the conditions 
$(U^1_x/(U^1)^2)_x=0$ and $(X_x/(U^1)^2)_x=0$. 
Denote the subgroup  of $G^\sim_0$, consisting of these transformations, by $\breve G^\sim_0$.    
The solutions of the system $S$ are constructed up to $\breve G^\sim_0$-equivalence. 

If $B_{xx}+2C_x\ne0$, the function~$\eta$ has the form $\eta=\eta^1(t,x)u$. Then $\eta^1=0$ and $C_x=0$ up to $\breve G^\sim_0$-equivalence. 
The interpretation of this solution is obvious. 
An associated ansatz for~$\mathcal L$ and the corresponding reduced equation are $u=\varphi(\omega)$, where $\omega=t$, and $\varphi_\omega=0$. 
The family of the associated invariant solutions of~$\mathcal L$ is formed by the constant functions. 

The condition  $B_{xx}+2C_x=0$ implies $C_{xx}=0$. Up to $\breve G^\sim_0$-equivalence we can assume that $B=C=0$.  
Then the system $S$ is reduced to the system $\eta_t=0$, $\eta_x+\eta\eta_u=0$. 
Its nonzero solutions are implicitly determined by the formula $u=x\eta+w(\eta)$, where $w=w(\eta)$ is an arbitrary smooth function of~$\eta$. 
An associated ansatz for the equation~$\mathcal L$ is found from the condition $u=xu_x+w(u_x)$ which is the Clairaut's equation with 
the implicit parameter~$t$. We choose the ansatz $u=\varphi(\omega)x+w(\varphi(\omega))$, where $\omega=t$. 
The corresponding reduced equation is $\varphi_\omega=0$, i.e., the associated invariant solutions of~$\mathcal L$ has the form 
$u=cx+w(c)$, where $c$ is an arbitrary constant. 

Let us emphasize that the obtained results have a compact form only due to the consideration up to $\breve G^\sim_0$-equivalence. 
\end{example}

Now we present a single example concerning the system ${\rm DE}_1(\mathcal L)$.  
In view of Corollary~\ref{CorollaryOnEquivOfRedOps1OfLPEs} we can assume without loss of generality that $g^3=0$ and, 
therefore, consider only the two first equations of the system ${\rm DE}_1(\mathcal L)$. 
The $G^\sim_1$-invariance of the equation $g^3=0$ additionally justifies this assumption. 

\begin{example}\label{ExampleOnRedOps1G20OfLPEs}
The constraint $g^2=0$ is invariant with respect to the transformations from the equivalence group $G^\sim_1$, in which $U^1=1$.
These transformations are presented by 
formulas~\eqref{EqGenFormOfTransOfInhomLPEs}, \eqref{EqTransOfCoeffsOfLPE} and~\eqref{EqTransOfCoeffsOfRedOps1ForLPEs}, 
where $U^1=1$ and $U^0=0$, and form the subgroup of $G^\sim_1$, denoted by $\breve G^\sim_1$.
Up to the $\breve G^\sim_1$-equivalence, the coefficient~$A$ can be assumed equal to 1.  
Imposing the conditions $g^2=g^3=0$ and $A=1$, we reduce ${\rm DE}_1(\mathcal L)$ to the system 
\begin{gather}\label{EqDE1-1forLPEsWithG20G30}
g^1_t-g^1_{xx}+2g^1g^1_x+(Bg^1)_x+B_t=0,
\\ \label{EqDE1-2forLPEsWithG20G30}
C_t+g^1C_x+2g^1_xC=0.
\end{gather} 

Equation~\eqref{EqDE1-1forLPEsWithG20G30} is linearized 
to the equation $w_t=w_{xx}+(Bw)_x$ by the generalization $g^1=-w_x/w-B$ of the Cole--Hopf substitution and then 
to the equation $v_t=v_{xx}+Bv_x$ by the subsequent substitution $w=v_x$. 
In the case $C=0$, the resulting substitution $g^1=-v_{xx}/v_x-B$ is the confinement of transformation~\eqref{EqTrans1BetweenInitialAndDetEqsForLPEs}
under the assumptions $v^3=0$, $v^2=1$ and $v^1=v$, where $v$ is a nonconstant solution of~$\mathcal L$. 

Equation~\eqref{EqDE1-2forLPEsWithG20G30} admits a double interpretation depending on a reading of the phrase 
``the equation~$\mathcal L$ possesses the reduction operator~$\p_t+g^1\p_x$''.
It can be considered either as an additional constraint for the function $g^1$ or an equation in the coefficient~$C$. 
Choosing the second alternative, we obtain $C=v_x{}^{\!2}\Phi(v)$ for some function~$\Phi=\Phi(v)$.

If $C=0$, equation~\eqref{EqDE1-2forLPEsWithG20G30} is an identity. 
Therefore, the equation~$\mathcal L$ admits any reduction operators of the form $\p_t-(v_{xx}/v_x+B)\p_x$, where 
$v=v(t,x)$ runs through the set of nonconstant solutions of~$\mathcal L$. 
The corresponding two-parametric solution family of~$\mathcal L$ is $u=c_1v(t,x)+c_2$.
\end{example}

\begin{note}
Since we do not initially specify values of the arbitrary elements and derive conditions on arbitrary elements depending 
on possessed reduction operators, the above examples have features of inverse problems of group analysis.  
Namely, we simultaneously describe both reduction operators with certain properties and values of arbitrary elements 
for which the corresponding equations admit such reduction operators. 
A similar inverse problem for generalized conditional symmetries of evolution equations is investigated in~\cite{Sergyeyev2002}. 
Due to possibilities on variation of arbitrary elements and application of equivalence transformations,
the problems of this kind essentially differ from the problem of finding reduction operators of a fixed equation. 
\end{note}

\section{Applications}\label{SectionOnReductionOpsOfLinearTransferEqs}

In Sections~\ref{SectionOnLinearizationOfDetEqsForredictionOpsOfLPEs}--\ref{SectionOnLieReductionsOfDetEqsOfReditionOpsOfLPEs}
``no-go'' statements of different kinds have been proved for the reduction operators of 
the equations from class~\eqref{EqGenLPE}. 
The term ``no-go'' has to be treated only as the impossibility of exhaustive solving of the problem 
or the inefficiency of finding Lie symmetries and Lie reductions of the determining equations.
At the same time, imposing additional (non-Lie) constraints on coefficients of reduction operators,
one can construct particular examples of reduction operators 
and then apply them to the construction of exact solutions of an initial equation.
Since the determining equations have more dependent or independent variables and, therefore, more degrees of freedom than the initial ones,
it is more convenient often to guess a simple solution or a simple ansatz
for the determining equations, which can give a parametric set of complicated solutions of the initial equations.
(A similar situation is for Lie symmetries of first-order ordinary differential equations.) 
It is the approach that was used, e.g., in~\cite{Gandarias2001} to construct exact solutions of a (nonlinear) fast diffusion equation 
with reduction operators having the zero coefficients of~$\p_t$. 
Earlier this approach was applied to the interesting subclass of class~\eqref{EqGenLPE}, 
consisting of the linear transfer equations of the general form
\begin{equation} \label{EqOfLinearTransfer}
u_t=u_{xx}+\frac{h(t)}{x}u_x.
\end{equation}
These equations arise, in particular, under symmetry reduction of the Navier--Stokes equations 
\cite{Fushchych&Popowych1994-1,Popovych1995,Popovych1997}. 
Investigation of reduction operators allowed us to construct series of multi-parametric solutions of equations~\eqref{EqOfLinearTransfer} 
and, as a result, wide solution families of the Navier--Stokes equations, parameterized by constants and functions of~$t$.

We consider class~\eqref{EqOfLinearTransfer} as an example showing possible ways 
of imposing nontrivial additional constraints to determining equations. 
This subclass is singled out from the whole class~\eqref{EqGenLPE} by the conditions on arbitrary elements 
$A=1$, $(xB)_x=0$ and $C=0$.

We fix an equation~$\mathcal L$ from class~\eqref{EqOfLinearTransfer}. 
The maximal Lie invariance algebra of~$\mathcal L$ is the algebra

\vspace{1ex}

1) $\langle u\p_u,\,f\p_u\rangle$ if $h\not=\const$;

\vspace{1ex}

2) $\langle\p_t,\,D,\,\Pi_h,\,u\p_u,\,f\p_u\rangle$ if $h=\const,$ $h\not\in\{0,2\}$; 

\vspace{1ex}

3) $\langle\p_t,\,D,\,\Pi_h,\, 2\p_x-hx^{-1}u\p_u,\, G_h,\,u\p_u,\,f\p_u\rangle$ if $\;h\in\{0,2\}$. 
\vspace{1ex}

\noindent
Here $D=2t\p_t+x\p_x$, $\Pi_h=4t^2\p_t+4tx\p_x-(x^2+2(1+h)t)u\p_u$, $G_h=2t\p_x-(x+htx^{-1})u\p_u$.
The function $f=f(t,x)$ runs through the set of solutions of~$\mathcal L$.
The case $h=2$ is reduced to the linear heat equation ($h=0$) by the transformation 
$\tilde t=t$, $\tilde x=x$ and $\tilde u=xu$, cf. Theorem~\ref{TheoremOnGroupClassificationOfLPEs}.
The inversection of the maximal Lie invariance algebras of equations from class~\eqref{EqOfLinearTransfer} coincides with 
$\langle u\p_u,\,\p_u\rangle$, i.e., the kernel Lie symmetry group of class~\eqref{EqOfLinearTransfer} consists of 
scalings and translations of~$u$.

It is easy to see that the equation~$\mathcal L$ possesses no nontrivial Lie symmetries and, therefore, no Lie reductions if $h\ne\const$. 
At the same time, non-Lie reduction operators can be found for an arbitrary value of~$h$.

Any reduction operator of~$\mathcal L$ with the nonzero coefficient of $\p_t$ is 
$G^\infty(\mathcal L)$-equivalent to an operator $\p_t+g^1\p_x+g^2u\p_u$,
where the functions $g^1=g^1(t,x)$ and $g^2=g^2(t,x)$ satisfy the two first equations 
of the corresponding determining system ${\rm DE}_1(\mathcal L)$. 
Following Example~\ref{ExampleOnRedOps1G20OfLPEs}, we impose the additional constraint $g^2=0$. 
Then the second equation of ${\rm DE}_1(\mathcal L)$ is identically satisfied. 
The first equation of ${\rm DE}_1(\mathcal L)$ is rewritten in the form
\[
(g^1+hx^{-1})_t=(g^1_x-g^1(g^1+hx^{-1}))_x.
\]
We put the left and right hand sides equal to~0. Then $g^1=\chi(x)-hx^{-1}$, $g^1_x-g^1(g^1+hx^{-1})=\psi(t)$.
The compatibility of these equations implies that $\chi=-x^{-1}$ and $\psi=0$, i.e., $g^1=-(h(t)+1)x^{-1}$ 
and the corresponding reduction operator is
\[
Q=\p_t-(h(t)+1)x^{-1}\p_x.
\]
As a result, the equation~$\mathcal L$ possesses the family of $Q$-invariant solutions 
\begin{equation} \label{Eq2ndOrderPolynomialSolutionOfEqOfLinearTransfer}
u=c_2\left(x^2+2\int(h(t)+1)dt\right)+c_1.
\end{equation} 
 
Each reduction operator of~$\mathcal L$ with the zero coefficient of $\p_t$ is equivalent to an operator $\p_x+\eta\p_u$,
where the coefficient $\eta=\eta(t,x,u)$ satisfies the corresponding determining equation ${\rm DE}_0(\mathcal L)$:
\begin{equation}\label{EqDet0ForRedOpsOfEqOfLinearTransfer}
\eta_t=\eta_{xx}+2\eta\eta_{xu}+\eta^2\eta_{uu}+h(x^{-1}\eta)_x.
\end{equation}
Suppose that the same operator $\p_x+\eta\p_u$ is a reduction operator of all equations from class~\eqref{EqOfLinearTransfer}, 
i.e., the function  $\eta$ is a solution of~\eqref{EqDet0ForRedOpsOfEqOfLinearTransfer} for any value of~$h$.
This demand leads to the additional constraint $(x^{-1}\eta)_x=0$ implying that $\eta=x\zeta(t,u)$. 
We substitute the expression for $\eta$ into~\eqref{EqDet0ForRedOpsOfEqOfLinearTransfer} and split with respect to~$x$.  
Integrating the obtained system $\zeta_{uu}=0$, $\zeta_t=2\zeta\zeta_u$, we construct all its solutions: 
\[
\zeta=-\frac{u+\mu}{2(t+\varkappa)} \quad\mbox{or}\quad \zeta=\nu,
\]
where $\mu$, $\varkappa$ and $\nu$ are arbitrary constants. 
In other words, the common reduction operators of equations from class~\eqref{EqOfLinearTransfer} are exhausted, 
up to equivalence with respect to the kernel Lie symmetry group (more precisely, up to translations of~$u$), by the operators of the form 
\[
G_\varkappa=(2t+\varkappa)\p_x-xu\p_u \quad\mbox{and}\quad \p_x+\nu\p_u.
\]
(It is obvious that there are no common reduction operators with nonzero coefficients of $\p_t$.) 
The constant~$\varkappa$ cannot be put equal to 0 similarly to the constant~$\mu$
since translations of~$t$ do not belong to the kernel Lie symmetry group of class~\eqref{EqOfLinearTransfer} 
and the classification up to the equivalence group of class~\eqref{EqOfLinearTransfer} in not convenient for the consideration. 
The operator~$G_\varkappa$ is represented as the linear combination $G+\varkappa\p_x$ 
of the Galilean operator $G=2t\p_x-xu\p_u$ and translational operator $\p_x$. 
The non-reduced form for the coefficient of~$\p_x$ in~$G_\varkappa$ is chosen to obtain this representation. 
For any equation~$\mathcal L$ from class~\eqref{EqOfLinearTransfer} the reduction operator $R=\p_x+\nu\p_u$ is $G^\infty(\mathcal L)$-equivalent 
to the operator $\p_x$ which is trivial since the arbitrary element~$C$ equals 0 in class~\eqref{EqOfLinearTransfer}.
Another formulation of above result is the following: 
Each equation from class~\eqref{EqOfLinearTransfer} is conditionally invariant with respect to arbitrary linear combinations 
of the Galilean operator $G$ and the translational operator $\p_x$.
The family of $G_\varkappa$-invariant solutions of an equation of the form~\eqref{EqOfLinearTransfer} consists of the functions 
\[
u=c_1\exp\left\{-\frac{x^2}{2(2t+\varkappa)}-\int\frac{h(t)+1}{2t+\varkappa}dt\right\}.
\] 
The corresponding family for the operator $\p_x+\nu\p_u$ has the form~\eqref{Eq2ndOrderPolynomialSolutionOfEqOfLinearTransfer} 
with $c_2=\nu$.

The constructed exact solutions are generalized to series of similar solutions 
\[
u=\sum_{k=0}^NT^k(t)x^{2k}, \quad
u=\sum_{k=0}^NS^k(t)\left(\frac{x}{2t+\varkappa}\right)^{2k}
\exp\left\{-\frac{x^2}{2(2t+\varkappa)}-\int\frac{h(t)+1}{2t+\varkappa}dt\right\}.
\]
The functions $T^k=T^k(t)$ and $S^k=S^k(t)$ respectively satisfy systems of ODEs 
\begin{gather*}
T^k_t=2(k+1)(h(t)+2k+1)T^{k+1}, \quad k=\overline{0,N-1}, \quad T^N_t=0,
\\
S^k_t=2(k+1)(h(t)+2k+1)(2t+\varkappa)^{-2}S^{k+1},\quad k=\overline{0,N-1}, \quad S^N_t=0,
\end{gather*}
which are easily integrated. 
These series of exact solutions also can be found using different techniques connected with reduction operators 
and their generalizations, in particular, via nonlocal transformations in class~\eqref{EqOfLinearTransfer}, 
associated with reduction operators~\cite{Fushchych&Popowych1994-1,Popovych1995}.

\section{Discussion}\label{SectionOnDiscussionOfReductionOpsOfLPEs}

The main result of the present paper is the chain of ``no-go'' statements on reduction operators 
of linear $(1+1)$-dimensional parabolic equations. 
These statements shows that application of conventional methods to solving of the determining equations 
for coefficients of such operators cannot lead to reduction operators giving new exact solutions of initial equations. 
In both the cases naturally arising under the consideration, 
the determining equations form well-determined systems whose solving is in fact equivalent to solving of 
the corresponding equations from class~\eqref{EqGenLPE}. 
All transformational and symmetry properties of the determining equations are induced by the corresponding properties of 
the initial equations. 
Reduction operators constructed via Lie reductions of the determining equations are also connected with 
Lie invariance properties of the initial equations. 
Nevertheless, it is demonstrated in Section~\ref{SectionOnReductionOpsOfLinearTransferEqs} that 
involvement of ingenious empiric approaches different from the Lie one can give 
reduction operators which are useful for the construction of non-Lie exact solutions of equations from class~\eqref{EqGenLPE}.

Techniques developed in this paper can be applied to the general class of $(1+1)$-dimensional evolution equations. 
We also plan to consider generalized reduction operators 
of linear $(1+1)$-dimensional parabolic equations, whose coefficients depend on derivatives of~$u$. 
An interesting subject related to this is the connection between (generalized) reduction operators and Darboux transformations. 
Here we give some hints on this connection.

Consider a fixed tuple of linearly independent functions~$(\psi^1,\dots,\psi^p)$ 
of $t$ and~$x$, and the linear independence are assumed over the ring of smooth functions of~$t$.
The \emph{Darboux transformation} constructed with the tuple~$(\psi^1,\dots,\psi^p)$ 
is denoted by ${\rm DT}[\psi^1,\dots,\psi^p]$ and is defined by the formula 
\cite{Matveev&Salle1991,Popovych&Kunzinger&Ivanova2007}
\[
\tilde u={\rm DT}[\psi^1,\dots,\psi^p](u)=\frac{W(\psi^1,\dots,\psi^p,u)}{W(\psi^1,\dots,\psi^p)}.
\]
Here $W(\varphi^1,\ldots,\varphi^s)$ denote the Wronskian of the functions $\varphi^1$, \ldots, $\varphi^s$
with respect to the variable~$x$, i.e., $W(\varphi^1,\ldots,\varphi^s)=\det(\p^{i-1}\varphi^j/\p x^{i-1})_{i,j=1}^{\;s}$.
The initial ($u$) and, therefore, obtained ($\tilde u$) functions also depend on $t$ and~$x$. 

The transformation ${\rm DT}[\psi^1,\dots,\psi^p]$ is represented as the action of a linear $p$-order differential operator 
with differentiations with respect to only~$x$, ${\rm DT}[\psi^1,\dots,\psi^p](u)={\rm DT}[\psi^1,\dots,\psi^p]u$. 
The operator will be denoted by the same symbol as the transformation  
and called the \emph{Darboux operator} associated with the tuple $(\psi^1,\dots,\psi^p)$. 
In the cases $p=1$ and $p=2$ the expressions of the Darboux operators respectively are 
\[
{\rm DT}[\psi^1]=\p_x-\frac{\psi_x}\psi,\quad 
{\rm DT}[\psi^1,\psi^2]=\p_{xx}-\frac{(W(\psi^1,\psi^2))_x}{W(\psi^1,\psi^2)}\p_x+\frac{W(\psi^1_x,\psi^2_x)}{W(\psi^1,\psi^2)}.
\]

If the functions~$\psi^1$, \dots, $\psi^p$ are linearly independent solutions of an equation~$\mathcal L$ from class~\eqref{EqGenLPE} then 
they are linearly independent over the ring of smooth functions of~$t$~\cite{Popovych&Ivanova2004ConsLaws,Popovych&Kunzinger&Ivanova2007}. 
The Darboux transformation ${\rm DT}[\psi^1,\dots,\psi^p]$ maps the equation~$\mathcal L$ to the equation~$\tilde{\mathcal L}$ 
also belonging to the class~\eqref{EqGenLPE} and 
having the following values of arbitrary elements~\cite{Matveev&Salle1991,Popovych&Kunzinger&Ivanova2007}
\[
\tilde A=A,\quad
\tilde B=B+pA_x,\quad
\tilde C=C+pB_x+\frac{p(p+1)}2A_{xx}+\frac{W_x}W A_x+2\left(\frac{W_x}W\right)_x A,\quad
\]
where the abbreviation $W=W(\psi^1,\dots,\psi^p)$ is used.

Suppose that a reduction operator~$Q$ of~$\mathcal L$ has the canonical form 
and is associated with a first-order linear differential operator~$\widetilde Q$ acting on functions of~$t$ and~$x$.
It means that either $Q=\p_t+g^1\p_x+g^2u\p_u$ if $Q\in\mathcal Q_1(\mathcal L)$ or $Q=\p_x+\eta^1u\p_u$ if $Q\in\mathcal Q_0(\mathcal L)$. 
(Here $g^1$, $g^2$ and $\eta^1$ are functions of~$t$ and~$x$.)
In the first case the operator $\widetilde Q=-\p_t-g^1\p_x+g^2$ equals the operator $-A\,{\rm DT}[v^1,v^2]$ 
on the solution set of the equation~$\mathcal L$,  
where the solutions $v^i=v^i(t,x)$, $i=1,2$, of~$\mathcal L$ are determined according to Corollary~\ref{CorollaryOnLinearizationOfDetEqs1ForRedOpsOfLPEs}.
In the second case the coefficient~$\eta^1$ admits the representation $\eta^1=\Psi_x/\Psi$, where $\Psi=\Psi(t,x)$ is a solution of~$\mathcal L$.
Therefore, $\widetilde Q=-{\rm DT}[\Psi]$. Finally, we have the following statement. 

\begin{proposition}
Let a reduction operator~$Q$ of an equation~$\mathcal L$ from class~\eqref{EqGenLPE} be associated, 
up to the equivalence relations of operators, 
with a first-order linear differential operator acting on functions of~$t$ and~$x$.  
Then it is equivalent to a Darboux operator constructed with one (resp. two) linearly independent solutions of this equations 
in the case of vanishing (resp. nonvanishing) coefficient of~$\p_t$.
\end{proposition}

The properties of single reduction operators of multi-dimensional equations essentially differ from that in 
the $(1+1)$-dimensional case. 
For example, all single reduction operators of $(1+n)$-dimensional linear heat equations are exhaustively classified in~\cite{Popovych&Korneva1998} 
for arbitrary~$n$ without addressing the general solution of this equation that annuls the possibility of ``no-go'' statements. 
At the same time, it is not the case for involutive families of reduction operators~\cite{Popovych1998,Vasilenko&Popovych1999}.

\subsection*{Acknowledgements}

The research was supported by the Austrian Science Fund (FWF), START-project Y237 and Lise Meitner project M923-N13.
The author is grateful to Vyacheslav Boyko and Michael Kunzinger  for useful discussions and interesting comments 
and also wish to thank the referees for their suggestions for the improvement of this paper.

\end{document}